\documentclass[final]{siamltex}

\usepackage[utf8]{inputenc}
\usepackage[T1]{fontenc}
\usepackage{amsmath,pifont}
\usepackage{amssymb}
\usepackage{color,pgf}
\usepackage{yfonts}
\usepackage{mathrsfs,textcomp}
\usepackage{enumerate}
\usepackage{latexsym}
\usepackage[english]{babel}
\usepackage{graphicx}
\usepackage{verbatim}
\newcommand{\ppotR}[3]
{

\begin{figure}\begin{center}
\includegraphics[width=#3truecm]{#1.eps}\\
\caption{#2}
\llabel{#1}
\end{center}
\end{figure}
\noindent$\!\!$
}

\newcommand{\immagine}[3][8]{
\begin{figure}[htb]
\begin{center}
\includegraphics[width=#1cm]{#2.eps}
\caption{#3}
\label{fig:#2}
\end{center}
\end{figure}}

\newcommand{\PTR}{PT\R^2}
\newcommand{\G}{G}

\renewcommand{\L}{\mathfrak{L}}

\newcommand{\p}{\mathfrak{p}}
\renewcommand{\P}{\mathfrak{P}}

\renewcommand{\b}[1]{{\bf #1}}

\newcommand{\Dh}{\Delta_H}
\renewcommand{\span}[1]{\mathrm{span}\Pg{#1}}

\newcommand{\distr}{{\blacktriangle}}

\newcommand{\g}{{\mathbf g}}
\newcommand{\metr}{\g}

\newcommand{\mfunz}[5]{
$#1 : \left\{ \begin{array}{ccl}
 #2 & \rightarrow & #3 \\
 #4 & \mapsto& #5   \end{array} \right.$}

\newcommand{\mmfunz}[5]{
$$ #1 : \left\{ \begin{array}{ccl}
 #2 & \rightarrow & #3 \\
 #4 & \mapsto& #5   \end{array} \right. $$}
\newcommand{\fz}[3]{#1:\, #2 \rightarrow #3}

\newcommand{\llabel}[1]{{\label{#1}}}

\newcommand{\ffoot}[1]{}
\newcommand{\fffoot}[1]{}

\renewcommand{\r}[1]{(\ref{#1})}

\newcommand{\bi}{\begin{itemize}}
\newcommand{\ei}{\end{itemize}}
\newcommand{\be}{\begin{enumerate}}
\newcommand{\ee}{\end{enumerate}}

\newcommand{\bd}{\begin{description}}
\newcommand{\ed}{\end{description}}
\renewcommand{\i}{\item}

\newcommand{\bqn}{\begin{eqnarray}}
\newcommand{\eqn}{\end{eqnarray}}
\newcommand{\eqnn}{\nonumber\end{eqnarray}}
\newcommand{\eqnl}[1]{\llabel{#1}\end{eqnarray}}

\newcommand{\nn}{\nonumber\\}

\newcommand{\ba}[1]{\begin{array}{#1}}
\newcommand{\ea}{\end{array}}

\newcommand{\R}{\mathbb{R}}
\newcommand{\C}{\mathbb{C}}

\newcommand{\Z}{\mathbb{Z}}
\newcommand{\nei}{neighborhood}

\newcommand{\bproof}{\begin{proof}}
\newcommand{\eproof}{\end{proof}}
\newtheorem{Theorem}{\bf Theorem}
\newtheorem{Lemma}[Theorem]{\bf Lemma}
\newtheorem{Corollary}[Theorem]{\bf Corollary}
\newtheorem{Definition}[Theorem]{\bf Definition}
\newtheorem{Proposition}[Theorem]{\bf Proposition}
\newtheorem{remark}[Theorem]{\bf Remark}

\newcommand{\bt}{\begin{Theorem}}
\newcommand{\et}{\end{Theorem}}
\newcommand{\bl}{\begin{Lemma}}
\newcommand{\el}{\end{Lemma}}
\newcommand{\bp}{\begin{Proposition}}
\newcommand{\ep}{\end{Proposition}}
\newcommand{\bc}{\begin{Corollary}}
\newcommand{\ec}{\end{Corollary}}
\newcommand{\bdeff}{\begin{Definition}}
\newcommand{\edeff}{\end{Definition}}
\newcommand{\brem}{\begin{remark}\rm}
\newcommand{\erem}{\end{remark}}
\newcommand{\auth}[1]{{\sc #1}}
\newcommand{\tit}[1]{{\rm #1}}

\newcommand{\jou}[1]{{\it #1}}
\newcommand{\vol}[1]{{\it #1}}
\newcommand{\pp}[1]{pp.~#1}

\newcommand{\lam}{\lambda}
\newcommand{\al}{\alpha}
\newcommand{\eps}{\varepsilon}
\newcommand{\ga}{\gamma}
\newcommand{\de}{\delta}

\renewcommand{\th}{\theta}
\newcommand{\Id}{\mathrm{Id}}
\newcommand{\Lie}{\mathrm{Lie}}
\newcommand{\Pt}[1]{\left( #1 \right)}
\newcommand{\Pg}[1]{\left\{ #1 \right\}}
\newcommand{\Pq}[1]{\left[ #1 \right] }
\newcommand{\Pa}[1]{\langle #1 \rangle}

\newcommand{\dom}{\mathscr{D}}
\newcommand{\con}{\mathit{C}}

\newcommand{\ce}{\mathrm{ce}}
\newcommand{\se}{\mathrm{se}}
\newcommand{\Rep}{\mathfrak{X}^\lam}

\newcommand{\PPP}{{\cal I}}
\newcommand{\OOO}{{\cal D}}
\newcommand{\cecco}[1]{{ #1}}

\begin{document}
\title{Anthropomorphic image reconstruction via hypoelliptic diffusion\thanks{This research has been supported  by the European Research Council, ERC
StG 2009 ``GeCoMethods", contract number 239748, by the ANR ``GCM", program ``Blanc--CSD"
project number NT09-504490, and by the DIGITEO project ``CONGEO".
}}

\author{Ugo Boscain\thanks{CMAP, \'Ecole Polytechnique CNRS, France,
 and Team GECO, INRIA Saclay, {\tt boscain@cmap.polytechnique.fr}} \and Jean Duplaix\thanks{Laboratoire LSIS, Universit\'e de Toulon, France, {\tt duplaix@univ-tln.fr}} \and  Jean-Paul Gauthier\thanks{Laboratoire LSIS, Universit\'e de Toulon, France  and Team GECO, INRIA Saclay, {\tt gauthier@univ-tln.fr}} 
\and Francesco Rossi\thanks{Aix-Marseille Univ, LSIS, 13013, Marseille, France, {\tt francesco.rossi@lsis.org}
}}


\maketitle
\begin{abstract}
In this paper we study  a model of geometry of vision 
due to  Petitot, Citti and Sarti.
One of the main features of  this  model is that the primary visual cortex V1 lifts an image from $\R^2$ to the bundle of directions of the plane.
Neurons are grouped into {\it orientation columns}, each of them  corresponding to a point of this  bundle. 

In this model  a corrupted image is reconstructed by minimizing the energy necessary for the activation of the orientation columns corresponding to regions in which the image is corrupted. The minimization process intrinsically defines an hypoelliptic heat equation on  the bundle of directions of the plane.

In the original model, directions are considered both with and without orientation, giving rise respectively to a problem on the group of rototranslations of the plane $SE(2)$ or on the projective tangent bundle of the plane $PT\R^2$.

We provide a mathematical proof of several important facts for this model.
We first prove that the model is mathematically consistent only if directions are considered without orientation. We then 
prove that the convolution of a $L^2(\R^2,\R)$ function (e.g. an image) with a 2-D Gaussian is 
generically a Morse function. This fact is important since the  lift of Morse functions to $PT\R^2$ is defined on a smooth manifold. 
We then provide the explicit expression of the hypoelliptic heat kernel on $PT\R^2$ in terms of Mathieu functions.

Finally, we present the main ideas of an algorithm which allows to perform 
image reconstruction on real non-academic images. A very interesting point is that this algorithm is massively parallelizable and needs no information on where the image is corrupted.

\end{abstract}


\section{Introduction}

 
In this paper we study a model of geometry of vision due to Petitot, Citti and Sarti. The main reference for the model is the paper \cite{citti-sarti}. Its first version can be found in \cite{petitot, petitot-tondut}.  This model was also studied by the authors of the present paper in \cite{suzdal}, by Hladky and Pauls \cite{hladky} and, independently, by Duits et al. in a series of papers mostly for contour completion \cite{duits-vA} and contour enhancement  \cite{duits-q1,duits-q2}. This model has been called the {\it pinwheel model} by Petitot himself, see \cite{pinwheel}. See also \cite{petitot-libro,sanguinetti} and references therein.

 To start with, assume  that a grey-level image  is represented by a function $\PPP\in L^2(\OOO,\R)$, where $\OOO$ is an open bounded domain of $\R^2$. The algorithm that we present here is based on three crucial ideas coming from neurophysiology: 

\be
\i   It is widely accepted that the retina approximately smoothes the images by making the convolution with a Gaussian function (see for instance \cite{langer,occhio,gatto} and references therein), equivalently solving a certain isotropic heat equation. Moreover, smoothing by the same technique is a widely used method  in image processing. Then, it is an interesting question in itself to understand generic properties of these smoothed images. Our first result (proved in Appendix A) is that, given $G(\sigma_x,\sigma_y)$ the two dimensional Gaussian centered in $(0,0)$ with standard deviations  $\sigma_x,\sigma_y>0$, then the smoothed image 
$$
f=\PPP\ast\G(\sigma_x,\sigma_y)\in L^2(\R^2,\R)\cap\con^\infty(\R^2,\R),
$$
is generically a Morse function (i.e. a smooth function having as critical points only non-degenerate maxima, minima and saddles). This has interesting consequences, as explained in the following.

\brem In several applications, the convolution is made with a Gaussian of small standard deviations. Equivalently, the smoothed image can be obtained as the solution of an isotropic heat equation with small final time.
\erem

\brem These results can be generalized to non-Gaussian filters and even to non-linear smoothing processes.   See for instance \cite{damon} for some of these generalizations. 
\erem

\i {\bf The primary visual cortex V1 lifts the image from $\R^2$ to the bundle of directions of the plane $PT\R^2$.}

In a simplified model\footnote{For example, in this model we do not take into account the fact that the continuous space of stimuli is implemented via a discrete set of neurons.} (see {\cite{citti-sarti}} and \cite[p. 79]{petitot-libro}), neurons of V1 are grouped into {\it orientation columns}, each of them being sensitive to visual stimuli at a given point $a$ of the retina and for a given direction $p$ on it. The retina is modeled by the real plane, i.e. each point is represented by $a\in\R^2$, while the directions at a given point are modeled by the projective line, i.e. $p\in P^1$. Hence, the primary visual cortex V1 is modeled by the so called {\it projective tangent bundle} $\PTR:=\R^2\times P^1$. From a neurological point of view, orientation columns are in turn grouped into {\it hypercolumns}, each of them being sensitive to stimuli at a given point $a$ with any direction. In the same hypercolumn, relative to a point $a$ of the plane, we also find neurons that are sensitive to other stimuli properties, like colors.  In this paper, we  focus only on directions and therefore  each hypercolumn is represented by a fiber $P^1$ of the bundle $\PTR$.  See Figure \ref{fig:f-hyper-bis}.

\immagine[13]{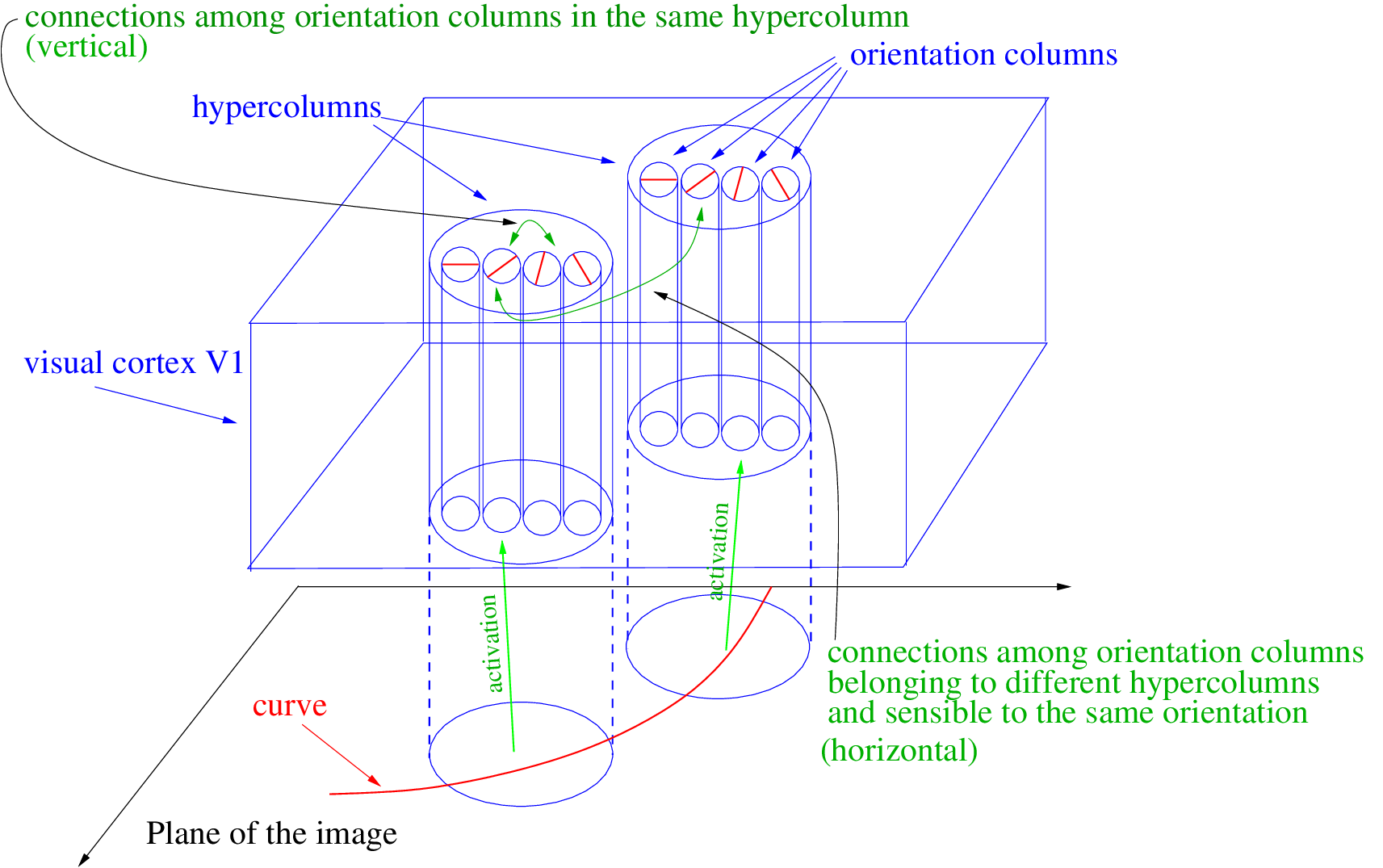}{A scheme of the primary visual cortex V1.}

This space has the topology of  $\R^2\times P^1$ (it is a trivial bundle) and its points are triples $(x,y,\th)$, where $(x,y)\in\R^2$, $\th\in\R/(\pi\Z)$.

 The smoothed image  $f:\R^2\to \R$ is {\it lifted} to a a function $\bar f$ defined as follows: 
$$\bar f(x,y,\th)=
\begin{cases}
f(x,y) &\mbox{\hspace{1mm} if $\th$ is the direction of the level set of $f$,}\\\\
0 & \mbox{\hspace{1mm}  otherwise.} 
\end{cases}
$$
It follows that  $\bar f$ has support  on a set $S_f\subset \PTR$.  The following fact constitutes our second result. 
 If $f$ is a Morse function (which happens generically due to the smoothing of the retina as explained above), then $S_f$ is  an embedded surface in  $\PTR$, see Proposition \ref{p-morse}.

\i If the image is corrupted or missing on a set $\Omega\subset\OOO$ (i.e. if $\PPP$ is defined on $\OOO\setminus\Omega$), then  the reconstruction  in $\Omega$ is made by minimizing a given cost. This cost represents the  energy that the  primary  visual cortex  should spend  in order to excite orientation columns which corresponds to points in $\Omega$ and hence that are not directly excited by the image. An orientation column is easily excited if it is close to another (already activated)  orientation column sensitive to a similar direction in a close position (i.e. if the two are close in $PT\R^2$).

When the image to be reconstructed is a curve, this gives rise to a sub-Riemannian problem (i.e. an optimal control problem  linear in the control and with quadratic cost) on $PT\R^2$, which we briefly discuss in Sections 
\ref{s-lift}, \ref{s-the-sub}, \ref{s-sachkov}. \footnote{In particular this sub-Riemannian structure has an underlying contact structure. To our knowledge, the first time in which the visual cortex was modeled as a contact structure was in \cite{hoffman}.}

When the image is not just a curve, the reconstruction is made by considering the diffusion process naturally associated with the sub-Riemannian problem on $PT\R^2$ (described by an hypoelliptic heat equation). Such a reconstruction  makes use of the function $\bar f$ as initial condition in a suitable way.   The reconstructed image is then obtained by projecting the result of the diffusion from $\PTR$ to $\R^2$.
\ee

In this paper we study this model, providing a mathematical proof of several key facts and adding certain important details with respect to its original version given in \cite{citti-sarti, sanguinetti}. 
 The main improvements are  the following:
\be[A) ] 
\i  As already mentioned, we start with any function $\PPP\in L^2(\OOO,\R)$ and we prove that after convolution with a Gaussian of standard deviations\footnote{We fix $\sigma_x=\sigma_y$ to guarantee invariance by rototranslations of the algorithm.}
 $\sigma_x=\sigma_y$, generically, we are left with a Morse function $f\in L^2(\R^2,\R)\cap\con^\infty(\R^2,\R)$ (see Appendix A). 
This smoothing process
is important to guarantee certain regularity of the domain of definition of the lifted function $\bar f$.

\i Our definition of the lift is suitable to all smooth functions, since we don't require conditions like nondegenerate gradient (as in \cite{citti-sarti}) or more complicated condition on the so called non-Legendrian solitary points (as in \cite[Thm 1.6]{hladky}).

\i Recall that $PT\R^2$ can be seen as the quotient of the group of rototranslations of the plane $SE(2)\simeq\R^2\times S^1$ by $\Z_2$, where the quotient is the identification of antipodal points in $S^1$. 
In the first version of this model \cite{citti-sarti} the image is lifted on $SE(2)$ (i.e. directions are considered with orientation), while in the second one \cite{sanguinetti} it is lifted on  $PT\R^2$ (i.e. directions are considered without orientation). 
The next contribution of our paper is to show that the problem of reconstruction of images for smooth functions
is well posed on $PT\R^2$ while it is not on $SE(2)$.
First, on $PT\R^2$ the lift is unique, while on $SE(2)$ it is not,  since level sets of the image are not oriented curves. 
Second, the problem on $SE(2)$ cannot be interpreted as a problem of reconstruction of contours (see Remarks \ref{r-existence}, \ref{r-se2} and \cite{suzdal}).
Third,   as proved in Proposition \ref{p-morse}, the domain $S_f$ of the lift of a Morse function $f$ is much more natural on $PT\R^2$ than on $SE(2)$. On $PT\R^2$ it is a manifold, while on $SE(2)$ it is a manifold with a boundary (for a continuous choice of the orientation of the level sets of $f$). The boundary appears on minima, maxima and saddles of $f$. In the diffusion process, starting with an initial condition which is concentrated on a manifold is much more natural than starting with an initial condition which is concentrated on a manifold with a boundary.

\i We show that the sub-Riemannian structure over $PT\R^2$ is not trivializable, which means that it cannot be specified by a single global orthonormal frame as in \cite{sanguinetti}.
For a detailed discussion of this issue see Remark \ref{r-2} and \cite{suzdal}.

\i   We give the expression of the hypoelliptic heat kernel over $PT\R^2$, while, previously, it  was known only on $SE(2)$ (see  \cite{nostro-kern} and  \cite{duits-vA,duits-q1}, where it  was found independently).
\i We provide an effective algorithm for image reconstruction that  looks unexpectedly efficient on real non-academic examples as shown in Section \ref{s-risultiagite}. Moreover, our algorithm  has the good feature to be massively parallelizable (see Section \ref{s-numerical}). This is just the materialization of the classical fact that the noncommutative Fourier transform disintegrates the regular representation over $SE(2)$. Moreover the algorithm does not need the information of where the original image is corrupted. 
\ee

Other numerical methods to compute hypoelliptic diffusion on $SE(2)$ for image processing have been developed. For instance:
group convolution methods (see \cite{libro-chir,duits-q1,franken-thesis}) 
finite differences \cite{citti-sarti, sanguinetti, franken-IJCV}\footnote{Notice that classical finite difference methods  ``hardly works'' to compute hypoelliptic diffusion. This is due to the 
diffusion at different scales on different directions as a consequence of the non-ellipticity of the diffusion operator.}  and finite elements methods. (See \cite{duits-q1,duits-q2}. These last works are related  to the noncommutative Fourier transform on $SE(2)$ and are extensions of the works by August \cite{august-thesis}.) Most of these works are about contour enhancement. 

\brem
 Notice that from the very beginning of the algorithm, we deal with the intensity of the image. Other related algorithms \cite[Sec. 3.3]{citti-sarti}, \cite{hladky} 
are instead composed of two reconstruction steps. After the lift of the image, these algorithms have to deal with a surface in $SE(2)$ or $\PTR$ with a hole, corresponding to the corrupted part. The first reconstruction step is thus to fill the hole as a surface, without considering the intensity of the image. The second reconstruction step is then to put the intensity on the reconstructed part. See Remarks \ref{r-difficulty} and  \ref{r-reiteriagite}.
\erem

%
%

 The results of our algorithm can be compared to the ones coming from psychological experiments. Moreover, they  can be useful to reconstruct the geometry of an image, as a preliminary step of  exemplar-based methods (see  \cite{masnou}).

Notice that an alternative technique of image processing (in particular for contour completion) based on physiological models of the visual cortex has been proposed by Mumford \cite{mumford}, then developed in \cite{duits-vA,duits-q1,duits-q2}.  In these models, contours (or level sets of an image) are considered with orientation and a
  non-isotropic diffusion is associated to an optimal control problem with drift having  elastica curves as solutions. We briefly compare the method presented in this paper with the one by Mumford in Section \ref{s-mumford}.




The structure of the paper is the following. In Section \ref{s-mat-alg} we present in detail the sub-Riemannian structure defined on $\PTR$.  We then define the corresponding hypoelliptic diffusion which is one of the main tools used in the algorithm of image reconstruction and we find explicitly the corresponding kernel on $\PTR$. At the end, we present in detail the mathematical algorithm.

Section \ref{s-numerical} is devoted to the discussion about  the numerical integration of the hypoelliptic evolution and to the presentation of samples.

%

 Appendix A is devoted to the detailed proof that, generically, the convolution of a $L^2$ function over a bounded domain $\OOO\subset\R^2$ with a Gaussian $G$ is a Morse function. In particular, we  prove  that the set of functions $\PPP\in L^2(\OOO)$ whose convolution with a Gaussian is a Morse function  in $L^2(\OOO)$ is residual (i.e. it is a countable intersection of open and dense sets). We then prove that the set of functions $\PPP\in L^2(\OOO)$ such that $\PPP\ast G$ restricted to a compact $K\subset\R^2$ is a Morse function is open and dense. Notice again that the proof can be adapted to any reasonable smoothing process, not only Gaussian.

\section{The mathematical model and the algorithm}
\label{s-mat-alg}

\subsection{Reconstruction of a curve}
\ffoot{qui si potrebbe pure pensare a descrivere mumford citti sarti etc come nel talk della rocchetta}
In this section we briefly describe an algorithm to reconstruct interrupted planar curves.
The main interest of this section is the definition of the sub-Riemannian structure over $PT\R^2$, from which we are going to define the sub-elliptic diffusion equation. {The main reference for this algorithm is \cite{citti-sarti}, where the lift of a planar curve was defined on $SE(2)$ rather than on $PT\R^2$}. 

Consider a smooth function $\fz{\ga_0}{[a,b]\cup[c,d]}{\R^2}$ (with $a<b<c<d$) representing a curve that is partially hidden or deleted in $(b,c)$. 
We assume that starting and ending points never coincide, i.e. $\ga_0(b)\neq\ga_0(c)$, and that initial and final velocities $\dot \gamma(b)$ and $\dot \gamma(c)$ are well defined and nonvanishing.

We want to find a curve $\fz{\ga}{\Pq{b,c}}{\R^2}$ that completes $\ga_0$ in the deleted part and that minimizes a cost depending both on the  length  and on the curvature $K_\gamma$ of $\gamma$.
Recall that 
$$
K_\gamma=\frac{\dot x \ddot  y-\dot  y\ddot  x}{(\dot  x^2+\dot  y^2)^{3/2}}$$

where $(x, y)$ are the components of $\gamma$.

The fact that $\ga$ completes $\ga_0$ means that $\ga(b)=\ga_0(b),\ \ga(c)=\ga_0(c)$. It is also reasonable to require that the directions of tangent vectors  coincide, i.e. $\dot{\ga}(b)\approx\dot{\ga_0}(b),\ \dot{\ga}(c)\approx\dot{\ga_0}(c)$ where 
\bqn
v_1\approx v_2 \mbox{ if it exists }\al\in\R\setminus\{0\} \mbox{ such that } v_1=\al\, v_2.
\eqnl{intro:ident-verso}
\brem
\label{r-existence}
Notice that we have required boundary conditions on initial and final directions without orientation. 
 The problem above  can also be formulated requiring 
boundary conditions with orientation, i.e. substituting in \r{intro:ident-verso} the condition $\al\in\R^+$.  However, this choice does not guarantee existence of minimizers for the cost we are interested with, see \cite{suzdal} and Remark \ref{r-2} below. 
An alternative model in which boundary conditions on directions are required with orientation is the one of Mumford. See Section \ref{s-mumford}.
\erem

In this paper we are interested in the minimization of the following cost, defined  for smooth curves $\gamma$ in $[b,c]$:
\bqn
J\Pq{\ga}=\int_b^c \sqrt{\| \dot{\ga}(t)\|^2 + \| \dot{\ga}(t)\|^2 K^2_\ga(t)}\, dt
\eqnl{eq-costocurva}

This cost is interesting for several reasons:
\bi
\i It depends on both  length and curvature of $\gamma$. It is small  for curves that are straight and short;

\i  It is invariant by rototranslation  (i.e. under the action of $SE(2)$) and by reparametrization of the curve, as should be any reasonable process of reconstruction of interrupted curves.

\i Minimizers for this cost do exist in the natural functional space in which this problem is formulated, without involving sophisticated functional spaces or curvatures that become measures. Indeed, in \cite{suzdal} we have proved:
\bp
For every $(x_b,y_b),(x_c,y_c)\in\R^2$ with $(x_b,y_b)\neq (x_c,y_c)$ and $v_b,\ v_c\in\R^2\backslash\Pg{0}$, the cost \r{eq-costocurva} has a minimizer over the set
\bqn
\dom&:=&\Pg{\ga\in\con^2([b,c],\R^2)\ \mbox{s.t.} 
\ba{l}
\sqrt{\| \dot{\ga}(t)\|^2(1+  K^2_\ga(t))}\in L^1([b,c],\R), \\
\ga(b)=(x_b,y_b),\ \ga(c)=(x_c,y_c),\\
\dot\ga(b)\approx v_b,\ \dot\ga(c)\approx v_c.
\ea
}.
\eqn
\ep

\brem
\label{r-2}
In \cite{yuri1,yuri2,yuri3}, it has been proved that minimizers for the cost \r{eq-costocurva} are analytic functions for which $\dot\gamma=0$ at most for two isolated points.  At these points the limit of  $\| \dot{\ga}(t)\| K_\ga(t)$ is well defined. They are cusp points, i.e. points at which $\dot\gamma$ becomes opposite. See Figure \ref{f-due-cusps}.

Notice that at cusp points the limit direction (regardless of orientation) is well defined. In \cite{suzdal} it is proved that if boundary conditions are required with orientation, then the cost \r{eq-costocurva}
has no minimum over the set $\dom$.
\erem

\ppotR{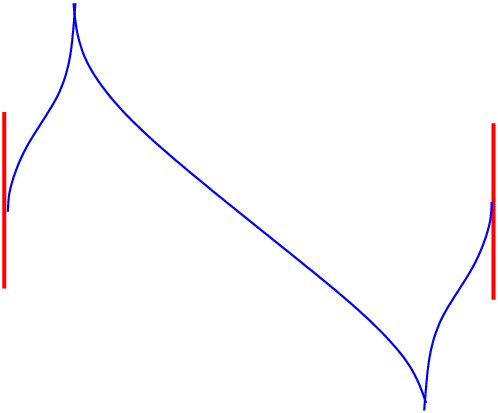}{A trajectory with two cusps.}{6}
\ei

However, the most interesting aspect from the modelling point of view is that this cost is a Riemannian length for lifts of planar curves 
over $PT\R^2$ (more precisely $J\Pq{\ga}$ is a sub-Riemannian length, see below). In the spirit of the model by Petitot, {Citti and Sarti}, this is the most natural distance that one can define on $\PTR$. Indeed, this distance takes into account the fact that two configurations $(x_1,y_1,\th_1)$ and $(x_2,y_2,\th_2)$ are close each other if they are both close in the planar coordinates $(x,y)$ and in the angle coordinate $\th$.

Apparently, this cost is a good model to describe the energy necessary to excite the orientation columns that are not directly excited by the image (since they correspond to the corrupted part of the image).
Indeed it is a standard fact in sub-Riemannian geometry (see Section \ref{s-www}) that the minimization of the cost $J\Pq{\ga}$ is equivalent to the minimization of  the energy-like cost
\bqn
E\Pq{\ga}=\int_b^c \left(\| \dot{\ga}(t)\|^2 + \| \dot{\ga}(t)\|^2 K^2_\ga(t)\, \right)dt.
\eqnn
The term $\| \dot{\ga}(t)\|^2$ models the energy necessary to activate {\it horizontal connections}, while  the term $\| \dot{\ga}(t)\|^2 K^2_\ga(t)$ models the energy necessary to activate {\it vertical connections}, see Figure \ref{fig:f-hyper-bis}. This is much more evident in the 
optimal control formulation of Section \ref{s-lift}, where $\| \dot{\ga}(t)\|^2$ is the control responsible for the ``straight movements on $\R^2$'' and $\| \dot{\ga}(t)\|^2 K^2_\ga(t)$ is the control corresponding for the ``rotational movements on $\R^2$''.  
Other models for these energies are of course possible,  but our choice appears to be the most natural since it provides a well-posed variational problem.

Finally,  a key consequence of this  choice of the cost is the following:  we have  a diffusion equation naturally associated with the  variational problem. This diffusion equation  can be used to reconstruct more complicated images than just curves. We  use this diffusion as the key tool for the reconstruction algorithm.

\brem One could argue that there is no reason to give the same weight to the length term $\|\dot \gamma \|$ and to the curvature term 
$\| \dot{\ga}(t)\|^2 K^2_\ga(t)$. However, if we define the cost 
\bqn
J_\beta[\ga]:=\int_b^c \sqrt{\| \dot{\ga}(t)\|^2 + \beta^2 \| \dot{\ga}(t)\|^2 K^2_\ga(t)}\, dt
\label{eq-beta}
\eqnn with a fixed $\beta\neq 0$ and if we consider an homothety $(x,y)\mapsto(\beta x,\beta y)$ and the corresponding transformation of a curve $\gamma=(x(t), y(t))$ to $\gamma_\beta=(\beta x(t),\beta y(t))$, then it is easy to prove that $J_\beta\Pq{\gamma_\beta}=\beta^2 J_1\Pq{\gamma}=\beta^2 J\Pq{\gamma}$. Therefore the problem of minimizing $J_\beta$ is equivalent to the minimization of $J$ with a suitable change of boundary conditions.

Although the mathematical problem is equivalent by changing $\beta$, this parameter will play a crucial role in the following, see Remark \ref{rem-beta}.
\erem

Another interesting feature is the uniqueness of this sub-Riemannian distance. Beside the possibility of adding a weight $\beta$ on the curvature term, that can be removed via an homothety, it is the unique sub-Riemannian distance for lift of planar curves on $PT\R^2$ that is invariant under the action of $SE(2)$. See Proposition \ref{p-!} below.
 
\subsection{Sub-Riemannian manifolds}\label{s-www}

In this section we recall some standard definitions of sub-Riemannian geometry, that we  use in the following. We start by recalling the definition of sub-Riemannian manifold.
\bdeff
A $(n,m)$-sub-Riemannian manifold is a triple $(M,\distr,{\mathbf g})$, 
where
\bi
\i $M$ is a connected smooth manifold of dimension $n$;
\i $\distr$ is a smooth distribution of constant rank $m< n$ satisfying the {\bf H\"ormander condition}, i.e. $\distr$ is a smooth map that associates to $q\in M$  a $m$-dim subspace $\distr(q)$ of $T_qM$ and $\forall~q\in M$ we have
\bqn
\hspace{-.5cm}\span{[X_1,[\ldots[X_{k-1},X_k]\ldots]](q)~|~X_i\in\mathrm{Vec}_H(M)}=T_qM
\eqnn 
where $\mathrm{Vec}_H(M)$ denotes the set of {\bf horizontal smooth vector fields} on $M$, i.e. $$\mathrm{Vec}_H(M)=\Pg{X\in\mathrm{Vec}(M)\ |\ X(q)\in\distr(q)~\ \forall~q\in M}.$$
\i $\g_q$ is a Riemannian metric on $\distr(q)$, that is smooth 
as function of $q$.
\ei
\edeff

A Lipschitz continuous curve $q(\cdot):[0,T]\to M$ is said to be {\it horizontal} if 
$\dot q(t)\in\distr(q(t))$ for almost every $t\in[0,T]$. Given an horizontal curve $q(\cdot):[0,T]\to M$, the {\it length of $q(\cdot)$} is
\bqn
l(q(\cdot))=\int_0^T \sqrt{ \g_{q(t)} (\dot q(t),\dot q(t))}~dt.
\eqnl{e-lunghezza}
The {\it distance} induced by the sub-Riemannian structure on $M$ is the 
function
\bqn
d(q_0,q_1)=&&\inf \{l(q(\cdot))\mid q(0)=q_0,q(T)=q_1,\ q(\cdot)\ \mathrm{horizontal}\}.
\eqnn 

The connectedness assumption for M and the H\"ormander condition guarantee the finiteness and the continuity of $d(\cdot,\cdot)$ with respect to the topology of $M$ (Chow's Theorem, see for instance \cite{agra-book}). The function $d(\cdot,\cdot)$ is called the {\it Carnot-Caratheodory distance} and gives to $M$ the structure of metric space (see \cite{bellaiche,gromov}).

It is a standard fact that $l(q(\cdot))$ is invariant under reparametrization of the curve $q(\cdot)$.
On one side, if an admissible curve $q(\cdot)$ minimizes the so-called {\it energy functional}
$$ E(q(\cdot))=\int_0^T {\g}_{q(t)}(\dot q(t),\dot q(t))~dt. $$
with  fixed $T$ (and initial and final fixed points), then $v=\sqrt{\g_{q(t)}(\dot q(t),\dot q(t))}$
is constant and $q(\cdot)$ is also a minimizer of $l(\cdot)$.
On the other side, a minimizer $q(\cdot)$ of $l(\cdot)$ such that  $v$ is constant is a minimizer of $E(\cdot)$ with $T=l(q(\cdot))/v$.

A {\it geodesic} for  the sub-Riemannian manifold  is a curve $q(\cdot):[0,T]\to M$ such that for each sufficiently small interval $[t_1,t_2]\subset [0,T]$, then $q(\cdot)_{|_{[t_1,t_2]}}$ is a minimizer of $E(\cdot)$.
A geodesic for which $\g_{q(t)}(\dot q(t),\dot q(t))$   is identically  equal to one is said to be {\it arclength parameterized.}

Locally, the pair $(\distr,{\mathbf g})$ can be specified by assigning a set of $m$ smooth vector fields spanning $\distr$, that are moreover orthonormal for ${\mathbf g}$, i.e.  
\bqn
\distr(q)=\span{X_1(q),\dots,X_m(q)},\qquad \metr_q(X_i(q),X_j(q))=\delta_{ij}.
\eqnl{trivializable}
Such a  set $\Pg{X_1,\ldots,X_m}$ is called a local {\it orthonormal frame} for the sub-Riemannian structure. When  $(\distr,{\mathbf g})$ can be defined by $m$  globally defined vector fields  as in \r{trivializable} we say that the sub-Riemannian manifold is {\it trivializable}.

Given a $(n,m)$-trivializable sub-Riemannian manifold, the problem of finding a curve minimizing the energy between two fixed points  $q_0,q_1\in M$ is
naturally formulated as the following optimal control problem 
\bqn
\dot q(t)&=&\sum_{i=1}^m u_i(t) X_i(q(t))\,,\label{sopra}\\
u_i(.)&\in& L^\infty([0,T],\R),\, \int_0^T \sum_{i=1}^m u_i^2(t)~dt\to\min,\nonumber\\ q(0)&=&q_0,~~~q(T)=q_1.
\eqnl{eq-op}
It is a standard fact that this optimal control problem is equivalent to the minimum time problem with controls $u_1,\ldots, u_m$ satisfying $u_1(t)^2+\ldots+u_m(t)^2\leq 1$ in $[0,T]$.
When the sub-Riemannian manifold is not trivializable, the equivalence with the optimal control problem \r{sopra}-\r{eq-op} is just local.


When the manifold is analytic and the orthonormal frame can be assigned by $m$ analytic vector fields, we say that the sub-Riemannian manifold is {\it analytic}.
In this paper we deal with an analytic sub-Riemannian manifold.

A sub-Riemannian manifold is said to be of 3D {\it contact} type  if $n=3$, $m=2$ and for every $q\in M$ we have span$\{\distr(q),[\distr,\distr](q)\}=T_qM$. This is the case that we study in this paper. For details, see \cite{agra-exp}.

\brem
As a consequence of the invariance by reparameterization of the cost \r{e-lunghezza}, it is equivalent to state the minimization problem in the space of  Lipschitz or absolutely continuous curves (i.e. for $u_i(\cdot)\in L^\infty([0,T],\R)$ or for $u_i(\cdot)\in L^1([0,T],\R)$.) See \cite[Lemma 1]{suzdal}.
\erem


\subsubsection{Left-invariant sub-Riemannian manifolds}
\label{ss-leftmanifold}

In this section we present a natural sub-Riemannian structure that can be defined on Lie groups. All along the paper, notations are adapted to group of matrices only. For general Lie groups, by $gv$ with $g$ in the Lie group $\G$ and $v$ in the Lie algebra $\L$, we mean $(L_g)_*(v)$ where $L_g$ is the left-translation on the group.

\bdeff 
Let $\G$ be a Lie group with Lie algebra $\L$ and $\P\subseteq\L$ a subspace of $\L$ satisfying the  Lie bracket generating condition $$\Lie~\P:=\span{[\p_1,[\p_2,\ldots,[\p_{n-1},\p_n]]]\ |\ \p_i\in\P}=\L.$$
Endow $\P$ with a positive definite quadratic form $\Pa{.,.}$. Define a sub-Riemannian structure on $\G$ as follows:
\bi
\i the distribution is the left-invariant distribution $$\distr(g):=g\P;$$
\i the quadratic form $\metr$ on $\distr$ is given by $$\metr_g(v_1,v_2):=\Pa{g^{-1}v_1,g^{-1}v_2}.$$
\ei
In this case we say that $(\G, \distr, \metr)$ is a left-invariant sub-Riemannian manifold.
\llabel{d-lieg-leftinv}
\edeff

In the following we define a left-invariant sub-Riemannian manifold by choosing a set of $m$ vectors $\Pg{\p_1,\ldots,\p_m}$ which form an orthonormal basis for the subspace $\P\subseteq\L$ with respect to the metric from  Definition \ref{d-lieg-leftinv}, i.e. $\P=\span{\p_1,\ldots,\p_m}$ and $\Pa{\p_i,\p_j}=\de_{ij}$. We thus have $$\distr(g)=g\P=\span{g\p_1,\ldots,g\p_m}$$ and $\g_g(g\p_i,g\p_j)=\de_{ij}$. Notice that every left-invariant sub-Riemannian manifold is trivializable.

\subsection{Lift of a curve on $PT\R^2$ and the sub-Riemannian problem}
\label{s-lift}
Consider a smooth planar curve $\gamma:[b,c]\to\R^2$.  This curve can be naturally lifted to a curve $\bar \gamma:[b,c]\to PT\R^2$ in the following way.
Let $(x(t),y(t))$ be the Euclidean coordinates of $\gamma(t)$. Then the coordinates of  $\bar\gamma(t)$ are  $(x(t),y(t),\th(t))$, where $\th(t)\in\R/(\pi\Z)$ is the direction of the vector $(x(t),y(t))$ measured with respect to the vector $(1,0)$. In other words,  
\bqn
\label{eq-arctan}
\th(t)=\arctan\Pt{\frac{\dot y(t)}{\dot x(t)}}\mbox{ mod }\pi.
\eqn

Of course we can extend by continuity the definition to points where $\dot\gamma(\bar t)=0$ but $\lim_{t\to \bar t}\th(t)$ is well defined. 
We assume
\bd
\i[{\bf [H]}] $\th:[b,c]\to\R/(\pi\Z)$ is absolutely continuous.
\ed
Notice that $\dot\theta=\|\dot\gamma\|K_\gamma$, hence hypothesis {\bf [H]} is equivalent to require that $\|\dot\gamma\|K_\gamma\in L^1([b,c],\R)$.

The requirement that a curve $(x(t),y(t),\th(t))$ satisfies the constraint \r{eq-arctan}
under {\bf [H]} can be slightly  generalized by  requiring  that 
$(x(t),y(t),\th(t))$ is an admissible trajectory of the control system on $PT\R^2$:
\bqn
&&\Pt{\ba{c}
\dot x\\
\dot y\\
\dot \th
\ea
}=u_1(t)\Pt{\ba{c}
\cos(\th)\\
\sin(\th)\\
0
\ea}+u_2(t)\Pt{\ba{c}
0\\
0\\
1
\ea}
\eqnl{eq-dyn}
with $u_1,u_2\in L^1([b,c],\R)$. Indeed each smooth trajectory $\gamma$ satisfying {\bf [H]} is an admissible trajectory of \r{eq-dyn}.

Since $u_1(t)^2=\|\dot \gamma(t)\|^2,$ $u_2(t)^2=\dot\th^2=\|\dot \gamma(t)\|^2 K_\gamma(t)^2$, we have
 \bqn
J\Pq{\gamma}=\int_b^c \sqrt{u_1(t)^2+ u_2(t)^2}\,dt
\eqnn 
Hence, the problem of minimizing the cost \r{eq-costocurva} on the set of curves $\dom$ is (slightly) generalized considering the optimal control problem 
(here $q(\cdot)=(x(\cdot),y(\cdot),\th(\cdot))$)
\bqn
&&\dot q=u_1(t)X_1(q)+u_2(t)X_2(q),\label{eq-dyn1}\\
&&X_1(q)=
\Pt{\ba{c}
\cos(\th)\\
\sin(\th)\\
0
\ea},~~X_2(q)=
\Pt{\ba{c}
0\\
0\\
1
\ea},\label{eq-dyn10}\\
&&l(q(\cdot))=\int_b^c \sqrt{u_1(t)^2+ u_2(t)^2}\,dt
\to\min,\label{eq-dyn2}\\
&&q(b)=(x_b,y_b,\th_b),\ q(c)=(x_c,y_c,\th_c),\label{eq-dyn3}\\
&&(x_b,y_b)\neq (x_c,y_c),\ u_1,u_2 \in L^1([b,c],\R).
\eqnl{eq-dyn4}

\brem
Notice that there are admissible trajectories $q(\cdot)=(x(\cdot),y(\cdot),\th(\cdot))$ of the control system \r{eq-dyn1} for which the condition $\th(\bar t)=\lim_{t\to\bar t}\arctan\Pt{\frac{\dot y( t)}{\dot x( t)}}$ is not verified
 (consider for instance the trajectory $x(t)=0$, $y(t)=0$, $\th(t)=t$) or such that 
$x(\cdot)$ or $y(\cdot)$ fail to be smooth.  However, it has been proved in \cite{suzdal} that minimizers of \r{eq-dyn1}-\r{eq-dyn4} are minimizers of \r{eq-costocurva} on the set $\dom$ and they are smooth.
\erem

\brem {\bf (non-trivializability)} A certain abuse of notation appears 
in  formulas \r{eq-dyn}, \r{eq-dyn2}, and \r{eq-dyn3}, as in \cite{sanguinetti}. Indeed the vector field 
$X_1$ is not well defined on $PT\R^2$. For instance, it takes two opposite values in $\th$ and $\th+\pi$, that are identified. A correct definition of the sub-Riemannian structure requires two charts:

\bi
\i Chart A: $\th\in]0+k\pi,\pi+ k\pi[$, $k\in\Z$, $x,y\in\R$.
 \bqn
&&\dot q=u_1^A(t)X_1^A(q)+u_2(t)X_2(q),~~X_1^A=
\Pt{\ba{c}
\cos(\th)\\
\sin(\th)\\
0
\ea},\nn
&&l(q(\cdot))=\int_b^c \sqrt{u_1^A(t)^2+ u_2(t)^2}\,dt,
\eqnn

\i Chart B: $\th\in]-\pi/2+k\pi,\pi/2+ k\pi[$, $k\in\Z$, $x,y\in\R$.
\bqn
&&\dot q=u_1^B(t)X_1^B(q)+u_2(t)X_2(q),~~X_1^B=
\Pt{\ba{c}
\cos(\th)\\
\sin(\th)\\
0
\ea},\nn
&&l(q(\cdot))=\int_b^c \sqrt{u_1^B(t)^2+ u_2(t)^2}\,dt,
\eqnn
\ei
One can check that the two charts are compatible and that this sub-Riemannian structure is non-trivializable, while  $PT\R^2$ is parallelizable.

Since  the formal expression of $X_1^A$ and   $X_1^B$ are the same, while they are defined on different  domains, one can proceed with a single  chart (however, one should  bear in mind that 
$u_1$ changes sign when passing from the chart A to the chart B in 
$\R\times\R\times]\pi/2,\pi[$). 
In the following, since we study a ``sum of squares'' hypoelliptic diffusion on this sub-Riemannian structure,  the problem disappears.
\erem

This sub-Riemannian manifold is of 3D contact type:  the distribution has dimension 2 over a three-dimensional manifold and $$\mathrm{span}\{X_1(q),\,X_2(q) ,\,[X_1,X_2](q)\}=T_{q}PT\R^2.$$ 
\ffoot{since $u_1$ and $u_2$ play the role of the component of $\gamma$ on the orthonormal base $F_1$ and $F_2$  and the cost \r{eq-dyn2} is simply the integral of the norm of $\dot \gamma$ in this base.
The name sub-Riemannian comes from the fact that the orthonormal base is formed by (noncommuting) vector fields in number less than the dimension of the space. The noncommutativity of the vector fields is crucial to guarantee the existence of a trajectory joining the boundary conditions. The term contact comes from the fact that span$\{F_1(x,y,\th),F_2(x,y,\th),[F_1,F_2](x,y,\th)\}=T_{(x,y,\th)}PT\R^2$.}

\subsection{The sub-Riemannian problem on $SE(2)$}
\label{s-the-sub}

It is convenient to lift the sub-Riemannian problem 
on $PT\R^2$ \r{eq-dyn1}-\r{eq-dyn4} on the group of rototranslation of the plane $SE(2)$, in order to take advantage of the group structure.  It is the group of matrices of the form 
\bqn
SE(2)=\Pg{\left(
\ba{ccc}
\cos(\th)&-\sin(\th)&x\\
\sin(\th)&\cos(\th)&y\\
0&0&1
\ea
\right)\ |\ 
\ba{l}
\th\in\R/(2\pi\Z),\\ 
x,y\in\R\ea}.
\eqnn
In the following we often denote an element of $SE(2)$ by $g=(x,y,\th)$.

A basis of the Lie algebra of $SE(2)$ is $\Pg{p_1,p_2,p_3}$, with
\bqn
p_1=\left(
\ba{ccc}
0&0 &1 \\
0&0 &0 \\
0 &0 &0 
\ea
\right), \hspace{.5cm}
p_2=\left(
\ba{ccc}
0&-1 &0 \\
1&0 &0 \\
0 &0 &0 
\ea
\right),\hspace{.5cm} 
p_3=\left(
\ba{ccc}
0&0 &0 \\
0&0 &1 \\
0 &0 &0 
\ea
\right).
\eqnn

We define a trivializable sub-Riemannian structure on $SE(2)$ as presented in Section \ref{ss-leftmanifold}: consider the two left-invariant vector fields $X_i(g)=g p_i$ with $i=1,2$ and set
\bqn
\distr(g)=\span{X_1(g),X_2(g)}~~~~~~\metr_g(X_i(g),X_j(g))=\de_{ij}.
\eqnn
In coordinates, the optimal control problem 
\bqn
&&\dot g \in\distr(g),~~l(g(.))=\int_b^c\sqrt{\metr_{g(t)}(\dot g,\dot g)}~dt\to\min,\label{eq-se1}\\
&&g(b)=(x^b,y^b,\th^b),~~g(c)=(x^c,y^c,\th^c),\\
&&(x^b,y^b)\neq (x^c,y^c),
\label{eq-se2}
\eqn
has  the form \r{eq-dyn1}-\r{eq-dyn4}, but  $\th \in \R/(2\pi\Z)$. Notice that the vector field 
$(\cos(\th),\sin(\th),0)$ is well defined on $SE(2)$.

\brem
\label{r-se2}
It is worth  mentioning that the problem \r{eq-se1}-\r{eq-se2} (i.e. the problem \r{eq-dyn1}-\r{eq-dyn4} with $\th \in \R /(2\pi\Z)$) cannot be interpreted as a problem of reconstruction of planar curves where initial and final positions and initial and final direction of velocities (with orientation) are fixed. For instance, consider the curve starting from $(x,y,\th)=(0,0,0)$ and corresponding to controls $u_1(t)=\pi/2-t$, $u_2(t)=1$. The corresponding trajectory in the $(x,y)$ plane  is $(-\cos (t)+\frac{1}{2} (\pi -2 t) \sin (t)+1,\pi
    \sin ^2\left(\frac{t}{2}\right)+t \cos (t)-\sin(t))$. Notice that this trajectory has a cusp at time $t=\pi/2$.
For $t\in[0,\pi/2[$ we have that $\th$ is the angle with respect to $(1,0)$ of the vector $(\dot x(t),\dot y(t))$, while for  $t\in]\pi/2,\pi]$, it is not. See Figure \ref{cuuusp}. 
\ppotR{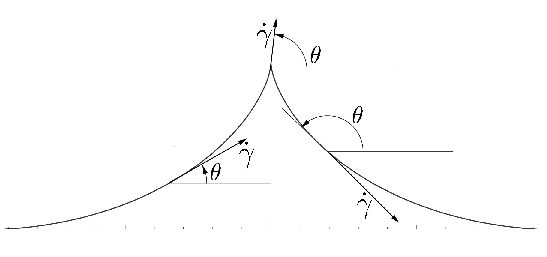}{A case in which $\th\in\R/(2\pi\Z)$ is not the direction (with orientation) of $\dot\ga$.}{7}
\erem

The control problem \r{eq-dyn1}-\r{eq-dyn4} defined on $\PTR$ is left-equivariant under the action of $SE(2)$. Indeed, topologically, $PT\R^2$ can be seen as the quotient of $SE(2)$ by $\Z_2$ (in other words $SE(2)$ is a double covering of $PT\R^2$). In coordinates, $(x,y,\th)\in PT\R^2$ corresponds to the two points $(x,y,\th),(x,y,\th+\pi)\in SE(2)$. Also, there is a natural transitive action of $SE(2)$ on $PTR^2$ given by
\bqn
\underset{\begin{smallmatrix} \in SE(2)
\end{smallmatrix}}
{\underbrace{\left(
\ba{ccc}
\cos(\th)&-\sin(\th)&x\\
\sin(\th)&\cos(\th)&y\\
0&0&1
\ea
\right)} }~
\underset{\in PT\R^2}{
\underbrace{\left(
\ba{c}
x'\\
y'\\
\th'
\ea
\right)}} =\underset{\in PT\R^2}{
\underbrace{\left(
\ba{c}
\cos(\th)x- \sin(\th)y+x'\\
\sin(\th)x+ \cos(\th)y+y'\\
\th'+\th
\ea
\right)}}
\nonumber
\eqn
where $\th'+\th$ is intended modulo $\pi$. The orthonormal frame for the sub-Riemannian structure  on $PT\R^2$ given by $X_1$ and $X_2$ in formula \r{eq-dyn1} is indeed left-equivariant under the action of $SE(2)$. 

In other words, given $(x,y,\th)\in\PTR$ such that $g\in SE(2)$ satisfies $(x,y,\th)=g(0,0,0)$, then 
\bqn
X_1(x,y,\th)=gp_1,~~ X_2(x,y,\th)=gp_2.
\eqn

\noindent
The following proposition can be checked directly.
{
\bp
\label{p-!}
Let $(PT\R^2,\distr,\metr)$ be a sub-Riemannian manifold and assume that it is left-equivariant under the natural action of $SE(2)$. This means  that if  $\{F_1,F_2\}$ is an ortnonormal frame for the sub-Riemannian structure  then 
\bqn
F_1(x,y,\th)=gF_1(0,0,0),~~ F_2(x,y,\th)=g F_2(0,0,0),
\eqn
where $g\in SE(2)$ is such that $(x,y,\th)=g(0,0,0)$.
Then, up to a change of coordinates and a rotation of the orthonormal frame, we have that 
\bqn
F_1(x,y,\th)=\left(  
\ba{c}
\cos(\th)\\ \sin(\th)\\0\ea\right).~~
F_2(x,y,\th)=
\left( 
 \ba{c}0\\0\\1/\beta
\ea
\right)
\label{Danilo-Restivo}
\eqn
for some $\beta>0$.
\ep
Notice that the problem of finding curves minimizing the length for the sub-Riemannian problem on $PT\R^2$ for which an orthonormal frame is given by \r{Danilo-Restivo}, is equivalent to the optimal control problem \r{eq-dyn1}, with the cost \r{eq-beta}.

}

\subsection{The Sachkov synthesis}
\label{s-sachkov}
The solution of the minimization problem \r{eq-dyn1}-\r{eq-dyn4} on $\PTR$, can be obtained from that of the problem on $SE(2)$ \r{eq-se1}-\r{eq-se2}. 
The latter has been studied by Yuri Sachkov in a series of papers \cite{yuri1,yuri2,yuri3}  (the first one in collaboration with I. Moiseev). 

 The author computed the optimal synthesis for the problem. More precisely he computed the geodesics starting from the identity and for each geodesic the cut time, i.e. the time where it loses optimality. Thanks to the group structure, optimal geodesics starting from other points are just translation of these ones. In Figure \ref{fig:f-yuri-mod} the cut locus of the Sachkov synthesis is shown.

\immagine[12]{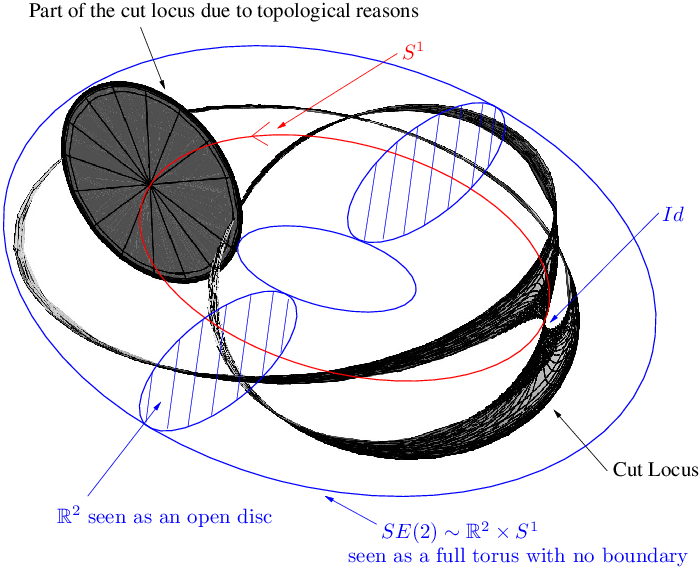}{The cut locus of the Sachkov synthesis, i.e. the set of points where geodesics lose optimality for the sub-Riemannian problem on $SE(2)$ (seen as the product of an open disc in $\R^2$ times $S^1$). Notice that the cut locus is adjacent to the starting point, as it always  occurs in sub-Riemannian geometry.}

The complete optimal synthesis and the description of the cut locus for the problem formulated on $PT\R^2$ has not been computed. However, as noticed by Sachkov, if we want to find the optimal trajectory joining   $(x,y,\th)$ to $(\bar x,\bar y,\bar \th)$ in $PT\R^2$, it is enough to find the shortest path among the four optimal trajectories joining the following points in $SE(2)$:
\bi
\i  $(x,y,\th)$ to $(\bar x,\bar y,\bar \th)$
\i $(x,y,\th+\pi)$ to $(\bar x,\bar y,\bar \th)$
\i $(x,y,\th)$ to $(\bar x,\bar y,\bar \th+\pi)$
\i  $(x,y,\th+\pi)$ to $(\bar x,\bar y,\bar \th+\pi)$
\ei

Moreover, Yuri Sachkov built a numerical algorithm for curve reconstruction on $PT\R^2$. In this paper, we will not go further on the subject of reconstruction of curves. For our purpose of image reconstruction,  the sub-Riemannian structure only is important, since it allows to define intrinsically a nonisotropic diffusion process.

\subsection{The hypoelliptic heat kernel}\label{s-HHK}

When the image is not just a curve, one cannot use the algorithm described above in which curves are reconstructed by solving a sub-Riemannian problem with fixed boundary conditions. Indeed, even if a corrupted image  
is thought as a set of interrupted curves (the level sets), it is unclear how to connect the different components of the level set among them (see Figure \ref{f-buco-frecce}).

Moreover, if the corrupted part contains the \nei\ of a maximum or minimum, then certain level sets are completely missing and cannot be reconstructed. 

\brem 
\label{r-difficulty}
 The difficulty of reconstructing a portion of an image containing a 
maximum or a minimum is also the main drawbacks of methods based on sub-Riemannian minimal surfaces. These algorithms (see \cite{citti-sarti,hladky,sanguinetti}) consider the boundary of the lift of the corrupted part as a closed curve $\gamma$ in the space $SE(2)$ or $\PTR$. They then ``fill the hole'' with the surface that has boundary $\gamma$ and minimizes the surface area computed with respect to the sub-Riemannian metric. As clearly explained in \cite{hladky}, this method can fail for 3 main reasons: the minimal surface does not exist (depending on the regularity of $\gamma$), it can be non-unique, or even can exist but its projection on $\R^2$ does not coincide with the corrupted part (either not covering the whole part or covering also a part of the non-corrupted image).

A second problem is that, even if the surface exists and it is computed, one has to choose how to diffuse the intensity of the image on the reconstructed surface. See  \cite[Def 7.3]{hladky} for the introduction of an ``interpolation'' function $f_t$ and a ``disambiguation'' function $F$. 
\erem

\ppotR{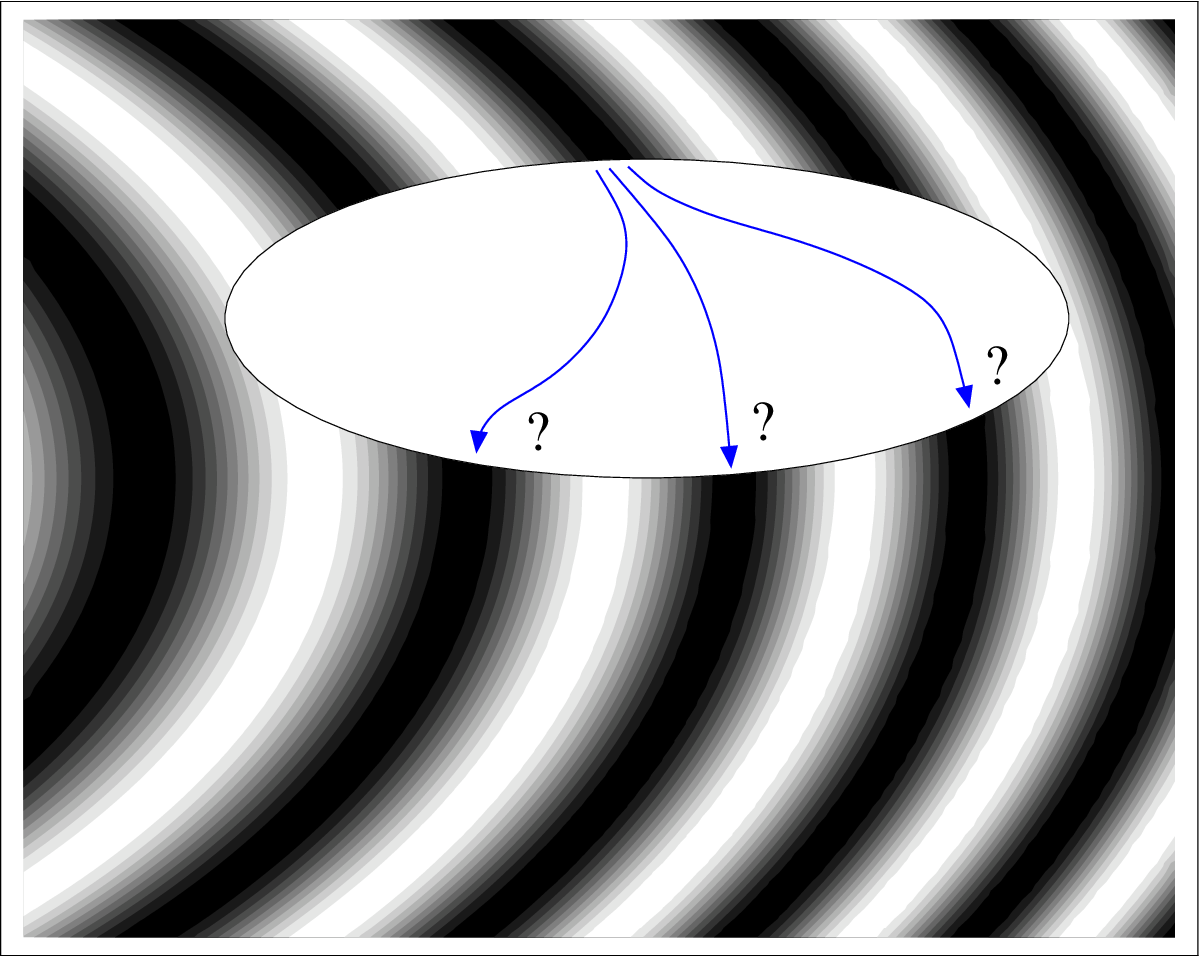}{The problem of connecting level sets}{8}

We then use the original image as the initial condition for the non-isotropic diffusion equation associated with the sub-Riemannian structure. { This idea was first presented in \cite{citti-sarti}.}

Roughly speaking,  first we consider all possible admissible paths by replacing the controls in equation \r{eq-dyn1} by independent Wiener processes. Then, we consider the  diffusion equation which describes the density of probability of finding the system in the point $(x,y,\th)$ at time $t$.

More precisely, let $\Pg{X_1,\ldots,X_m}$ be an orthonormal frame of a sub-Riemannian manifold and consider the 
stochastic differential equation $$dq_t =\sum_{i=1}^m X_i(q)_t\, dw^i_t,$$
where the $w^i$ are independent Wiener processes. It is a standard  result that, due to the H\"ormander condition, the stochastic  process admits a probability density $\phi(q,t)$, that satisfies the Fokker-Planck-like diffusion PDE\footnote{We emphasize here the fact that the PDE is not a stochastic one.}  
\bqn
\label{eq-heatgeneral}
\partial_t \phi(q,t)=\sum_{i=1}^m X_i^2 \phi(q,t).
\eqn 
 For more details, see e.g. \cite{meyer,pardoux}.

Roughly speaking,  the relation among the sub-Riemannian geodesics and the corresponding diffusion equation is the following: for small time the diffusion occurs mainly along optimal geodesics. 

For instance, a  result due to Leandre  \cite{leandremaj,leandremin}  states that if $p_t(q_1,q_2)$ is the heat kernel associated to \r{eq-heatgeneral} then 
for $t\to0$ we have that  $-t\log p_t(q_1,q_2)\to d(q_1,q_2)^2/4$, where  $d(\cdot,\cdot)$ is the  Carnot-Caratheodory distance.  For other results in this direction see \cite{benArous, jerison} and reference therein.

\smallskip
In our case, this diffusion equation is
\bqn
\partial_t \phi(x,y,\th,t)&=&\Delta_{H}\phi(x,y,\th,t)
\eqnl{eq-calore}
where
\bqn
\Delta_{H}&=&(X_1)^2+(X_2)^2=(\cos(\th)\partial_x+\sin(\th)\partial_y)^2+\partial_\th^2
\eqnn
Since at each point $(x,y,\th)$ we have $$\span{X_1, X_2, \Pq{X_1,X_2}}=T_{(x_1,x_2,\th)}PT\R^2,$$
 H\"ormander theorem \cite{horm} implies that  the operator $\Delta_{H}$ is hypoelliptic. 

The diffusion described by the equation \r{eq-calore} is highly non
isotropic. Indeed one can estimate the hypoelliptic heat kernel in terms of the sub-Riemannian distance (see for instance \cite{nostro-kern}), that is highly non isotropic as a consequence of the ball-box theorem (see for instance \cite{bellaiche}).

\brem
Notice that the sub-elliptic diffusion equation corresponding to the sub-Riemannian structure \r{eq-se1}-\r{eq-se2} on $SE(2)$, has the same form \r{eq-calore}. The only difference is that on $SE(2)$, $\th\in\R/(2\pi\Z)$. 
\erem
 
 \brem
In \cite{nostro-kern} it has been proved that the Laplacian  $\Delta_H$ is intrinsic on $SE(2)$, meaning  that it does not depend on the choice of the orthonormal frame for the sub-Riemannian structure. One can easily prove that this is the case also  for $\Delta_H$ on  $PT\R^2$.
 \erem
 
 \subsection{The hypoelliptic heat kernel on $SE(2)$}
The hypoelliptic heat kernel for the equation \r{eq-calore} on $SE(2)$ was  computed in \cite{nostro-kern, duits-vA}. More precisely, thanks to the left-invariance of $X_1$ and $X_2$, the equation \r{eq-calore} admits a a right-convolution kernel $p_t(.)$, i.e. there exists $p_t$ such that
\bqn
 e^{t \Dh } \phi_0(g)=\phi_0\ast p_t (g)=\int_G  \phi_0(h)p_t(h^{-1}g) \mu(h)
\eqnl{e:sol-SE(2)}
is the solution for $t>0$ of \r{eq-calore} with initial condition $\phi(0,g)=\phi_0(g)\in L^1(SE(2),\R)$ with respect to the Haar measure $\mu$.

We have computed $p_t$ on $SE(2)$ in \cite{nostro-kern}:
\bqn
\llabel{eq-M2-heat-exp}
p_t(g)&=&\int_0^{+\infty}\lam\left(
\sum_{n=0}^{+\infty} e^{a_n^\lam t} <\ce_n (\th,\frac{\lam^2}{4}),\Rep(g)\ce_n (\th,\frac{\lam^2}{4})>+\right.\nn
&&+\left.\sum_{n=1}^{+\infty} e^{b_n^\lam t} <\se_n (\th,\frac{\lam^2}{4}),\Rep(g)\se_n (\th,\frac{\lam^2}{4})>\right)\ d\lam.
\eqn
Here $\lam$ indexes the unitary irreducible representations of the group and 
\bqn
&&\Rep(g): L^2(S^1,\C)\to L^2(S^1,\C),\nn
&&\Rep(g)\psi(\al)=e^{i\lam(x\cos(\al)-y\sin(\al))}\psi(\al+\th)
\eqnn
 is the representation of the group element $g=(x,y,\th)$ on $L^2(S^1,\C)$.

 The functions se$_n$ and ce$_n$ are the $2\pi$-periodic Mathieu cosines and sines, and $<\phi_1,\phi_2>:=\int_{S^1}\phi_1(\al)\overline{\phi_2(\al)}\,d\al.$ The eigenvalues of the hypoelliptic Laplacian are $a_n^\lam:=-\frac{\lam^2}{4}-a_n\Pt{\frac{\lam^2}{4}}$ and $b_n^\lam:=-\frac{\lam^2}{4}-b_n\Pt{\frac{\lam^2}{4}}$, where $a_n$ and $b_n$ are characteristic values for the Mathieu equation. For details about Mathieu functions see for instance \cite[Chapter 20]{mathieu}. 

Since the operator $\partial_t-\Dh$ is hypoelliptic, then the kernel is a $\con^\infty$ function of $(t,g)\in \R^+\times G$. 
Notice that $p_t(g)=e^{t \Dh }\de_\Id(g)$.\\

The kernel \r{eq-M2-heat-exp} has been obtained by using the generalized Fourier transform. Once again, we refer to \cite{nostro-kern} for a detailed description of the generalized Fourier transform and the method to compute the kernel.

\subsection{The hypoelliptic heat kernel on $PT\R^2$}
$SE(2)$ is a double covering of $PT \R^2$. To a point $(x,y,\th)\in PT\R^2$  correspond the two points $(x,y,\th)$ and  $(x,y,\th+\pi)$ in $SE(2)$.  From the next proposition it follows that we can interpret the hypoelliptic heat equation on $PT\R^2$ as the hypoelliptic heat equation on $SE(2)$ with a symmetric initial condition. It permits also to compute the heat kernel on $PT\R^2$ starting from the one on $SE(2)$.

\bp
Let $\phi_0\in L^1(SE(2),\R)$ and assume that $\phi_0(x,y,\th)=\phi_0(x,y,\th+\pi)$ a.e. Then the solution at time $t$ of the hypoelliptic heat equation \r{eq-calore} on $SE(2)$,  having $\phi_0$ as initial condition at time zero, satisfies
\bqn\phi(t,x,y,\th)=\phi(t,x,y,\th+\pi).\eqnl{e:inizio-simm}
Moreover if  $\phi_0\in L^1(PT\R^2,\R)$, then the solution at time $t$ of the hypoelliptic heat equation on $PT\R^2$ \r{eq-calore} having $\phi_0$ as initial condition at time zero is given by
\bqn
\phi(t,x,y,\th)&=&\int_{PT\R^2}  \phi_0(\bar x,\bar y,\bar\th) P_t(x,y,\th, \bar x,\bar y,\bar\th)\, d\bar x\, d\bar y\, d\bar \th
\eqnl{e:kern-PTR-esp}
where
\bqn
P_t(x,y,\th,\bar x,\bar y,\bar\th)&:=&   p_t((\bar x,\bar y,\bar\th)^{-1}\circ(x,y,\th)) +  p_t((\bar x,\bar y,\bar\th)^{-1}\circ(x,y,\th+\pi)). 
\label{eq-kernelloni} 
\eqn
In the right hand side of equation \r{eq-kernelloni}, the group operations are intended in $SE(2)$.
\ep

\bproof
Define the element $\Pi=(0,0,\pi)\in SE(2)$ and observe the following properties:
\bi
\i $\Pi$ is idempotent.
\i Property \r{e:inizio-simm} reads as $\phi_0(g \Pi)=\phi_0(g)$.
\i The kernel $p_t(g)$ satisfies $p_t(\Pi g)=p_t(g \Pi)$. Indeed, call $g=(x,y,\th)$ and observe that, given a real function $\psi(\al)$, we have
\bqn\Rep\Pt{\Pi\circ g} \psi(\al)&=&\Rep\Pt{(-x,-y,\th)} \psi(\al)=\nn
&=&\overline{\Rep\Pt{(x,y,\th+\pi)} \psi(\al)}=
\overline{\Rep\Pt{g\circ \Pi} \psi(\al)}. 
\eqnn
\noindent
Recalling the explicit expression of $p_t$ given in \r{eq-M2-heat-exp}, we have $p_t(\Pi g)=\overline{p_t(g \Pi)}$. But $p_t$ is real, hence the equality follows.
\ei

We compute now $\phi(t,g \Pi)$ in $SE(2)$ with $\phi_0$ satisfying \r{e:inizio-simm} and we prove that $\phi(t,g \Pi)=\phi(t,g)$. Indeed,
\bqn 
\phi(t,g \Pi)&=&\int_G \phi_0(h) p_t(h^{-1} g \Pi)\,dh=
\int_G \phi_0(l \Pi) p_t(\Pi^{-1}l^{-1} g \Pi)\,d(l \Pi)=\nn
&=&\int_G \phi_0(l) p_t(\Pi \,\Pi^{-1} l^{-1} g) \,dl=\phi(t,g).
\eqnn

We now prove the expression  \r{e:kern-PTR-esp} for $\phi (t,\Pq{g})\in L^1(PT\R^2,\R)$ for initial data $\phi_0(\Pq{g})$, where $\Pq{g}$ is an element of $\PTR$, the class containing $g$ and $g\Pi$ in $SE(2)$. Consider the function $\psi_0(g)\in L^1(SE(2),\R)$ defined by $\psi_0(g)=\phi(\Pq{g})$, that clearly satisfies \r{e:inizio-simm}. Consider the unique solution $\psi(t,g)$ of the hypoelliptic equation \r{eq-calore}, that is given by $\psi(t,g)=\psi_0\ast p_t(g)$. Since $\psi(t,g)=\psi(t,g\Pi)$, the function $\phi(t,\Pq{g}):=\psi(t,g)$ is well defined.

It remains to show that $\phi$ defined above is the solution of \r{eq-calore} on $\PTR$. Indeed $\partial_t \phi=\partial_t \psi =\Dh \psi.$ Since the vector fields defining $\Dh$ both on $SE(2)$ and $\PTR$ coincide, then the differential operators $\Dh$ defined on $SE(2)$ and $\PTR$ coincide, hence $\Dh \psi=\Dh \phi$. Thus $\phi$ satisfies \r{eq-calore}. Since $\phi(0,\Pq{g})=\phi_0(\Pq{g})$, then $\phi$ is the (unique) solution.

The explicit expression \r{e:kern-PTR-esp} is a direct consequence of the definition $\phi(t,\Pq{g}):=\psi(t,g)$ and of the explicit expression of $\psi$ given in \r{e:sol-SE(2)}. Indeed,
\bqn  \phi(t,\Pq{g})&=&\psi(t,g)=\int_{SE(2)}  \psi_0(h)p_t(h^{-1}g) dh=
\int_{\R^2}\int_0^{2\pi}  \psi_0(h)p_t(h^{-1}g)\,dh=\nn
&=&\int_{\R^2}\int_0^{\pi}  \psi_0(h)p_t(h^{-1}g)+\psi_0(h\Pi)p_t((h \Pi)^{-1}g) \, dh=\nn
&=&\int_{\PTR} \phi_0(h)\Pt{p_t(h^{-1}g)+p_t(h^{-1}g\Pi)}\, dh.
\eqnn
The expression \r{e:kern-PTR-esp} is recovered by writing $g=(x,y,\th),\,h=(\bar x,\bar y,\bar\th)$ and recalling that $g\Pi=(x,y,\th+\pi)$.
\eproof

\subsubsection{Oriented vs. non-oriented approach}
\label{s-mumford}
One of the key points of the algorithm presented in this paper is that directions are considered without orientation. As mentioned above,  this choice is forced by well-posedness arguments when using the sub-Riemannian cost.
Other approaches which consider directions with orientation are possible, but with a different cost.

The most celebrated is the one due to Mumford \cite{mumford}, which in control language reads:
\bqn
&&\dot q=
\left(
\ba{c}
\cos(\th)\\
\sin(\th)\\
0
\ea
\right)+u(t)\left(
\ba{c}
0\\
0\\
1
\ea
\right),~~~\int_b^c (1+ \beta^2 u(t)^2)\,dt=\int_b^c (1+ \beta^2 K^2_{\gamma(t)})\,dt
\to\min,\nonumber\\
&&q(b)=(x_b,y_b,\th_b),\ q(c)=(x_c,y_c,\th_c),\nonumber
\eqnn
Here $(x,y,\th)\in\R^2\times S^1$. Notice that this variational problem has the same form as the one treated in this paper (i.e.  \r{eq-dyn1}-\r{eq-dyn10} for the energy functional $\int_b^c (u_1(t)^2+ \beta^2 u_2(t)^2)\,dt$), but forcing $u_1$ to be one and taking $\theta\in S^1$ instead of $P^1$. In this way trajectories are parametrized by arclength  for the Euclidean metric on $\R^2$ and not for the sub-Riemannian one and, as a consequence, there are no cusps. Notice, however, that since the energy functional is not invariant by reparametrization, geodesics have a different expression. (These geodesics are ``elastica'' curves, see also \cite{yuri-elasta1,yuri-elasta2}.)

 With  a procedure similar to the one described in Section \ref{s-HHK},  one can 
 naturally
 associate 
  a diffusion equation
 to this model. 
  Then, one gets a diffusion equation with drift, namely
$$
\partial_t \phi(x,y,\th,t)=\left(X_1+(X_2)^2\right)\phi(x,y,\th,t)=\left((\cos(\th)\partial_x+\sin(\th)\partial_y+\partial_\th^2\right) \phi(x,y,\th,t).
$$
{The level set of the image can be oriented for instance on the left (or on the right) of the gradient of the initial condition. This choice is well defined when the initial condition is a Morse function. In practice, people consider both  diffusions with positive and negative drift to have a more ``symmetric'' impainting.  This approach was followed in \cite{duits-q1,duits-q2} for contour enhancement.   }
 
Apparently, in the community, some researchers prefer  Mumford's model, while others prefer the Petitot-Citti-Sarti model presented in this paper. Mumford's model has the advantage of not producing cusps (which are not observed in psychological experiments, see \cite{petitot-libro}), while the  model presented in this paper has the advantage of treating horizontal and vertical connections at the same level  and allows a more natural lift of the image.



\subsection{The mathematical algorithm}
In this section we describe the main steps of the mathematical algorithm for image reconstruction. In the next section we give some guidelines for numerical implementation.

\begin{figure}[htb]
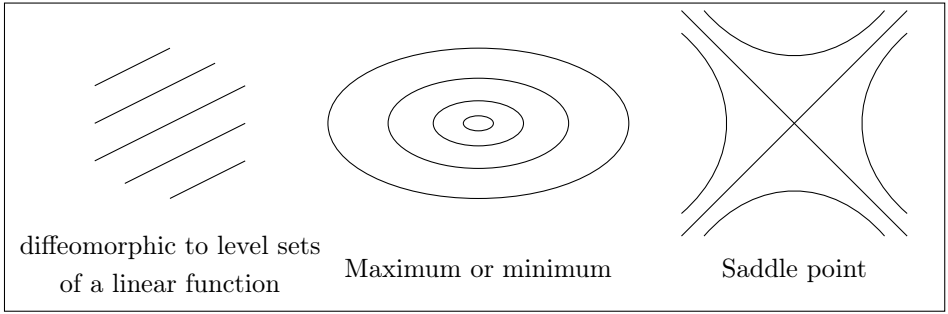

\begin{center} \begin{pgfpictureboxed}{0cm}{0cm}{12.5cm}{4.1cm}
\pgfline{\pgfxy(1.4+.2,1.2+.5)}{\pgfxy(3+.2,2+.5)}
\pgfline{\pgfxy(2+.2,1+.5)}{\pgfxy(3+.2,1.5+.5)}
\pgfline{\pgfxy(1+.2,1.5+.5)}{\pgfxy(3+.2,2.5+.5)}
\pgfline{\pgfxy(1+.2,2+.5)}{\pgfxy(2.6+.2,2.8+.5)}
\pgfline{\pgfxy(1+.2,2.5+.5)}{\pgfxy(2+.2,3+.5)}
\pgfputat{\pgfxy(2+.2,.8+.2)}{\pgfbox[center,top]{diffeomorphic to level sets}}
\pgfputat{\pgfxy(2+.2,.3+.2)}{\pgfbox[center,top]{of a linear function}}

\pgfellipse[stroke]{\pgfxy(6+.3,2+.5)}{\pgfxy(0,1)}{\pgfxy(2,0)}
\pgfellipse[stroke]{\pgfxy(6+.3,2+.5)}{\pgfxy(0,.6)}{\pgfxy(1.2,0)}
\pgfellipse[stroke]{\pgfxy(6+.3,2+.5)}{\pgfxy(0,.3)}{\pgfxy(.6,0)}
\pgfellipse[stroke]{\pgfxy(6+.3,2+.5)}{\pgfxy(0,.1)}{\pgfxy(.2,0)}
\pgfputat{\pgfxy(6.3,.7)}{\pgfbox[center,top]{Maximum or minimum}}

\pgfline{\pgfxy(9,1)}{\pgfxy(12,4)}
\pgfline{\pgfxy(9,4)}{\pgfxy(12,1)}
\pgfxycurve(9,1.3)(9.8,2)(9.8,3)(9,3.7)
\pgfxycurve(12,1.3)(11.2,2)(11.2,3)(12,3.7)
\pgfxycurve(9.3,1)(10,1.8)(11,1.8)(11.7,1)
\pgfxycurve(9.3,4)(10,3.2)(11,3.2)(11.7,4)
\pgfputat{\pgfxy(10.5,.7)}{\pgfbox[center,top]{Saddle point}}

\end{pgfpictureboxed}
\end{center}
\caption{Level sets of a Morse function}
\label{morse}
\end{figure}

\vskip 1cm\noindent
{\bf STEP 1: Smoothing of $\PPP_c$}
Assume that the grey level of a corrupted image is described by a function $\PPP_c:\OOO_c:=\OOO^2\setminus \Omega \to[0,\infty[$.
The set $\Omega$ represents the region where the image is corrupted.  The subscript ``$_c$'' means ``corrupted''. 
After making the convolution with a Gaussian of standard deviations $\sigma_x=\sigma_y>0$\footnote{$\PPP_c$ is considered to be zero outside $\OOO_c$. Moreover we assume $\sigma_x=\sigma_y$ to guarantee invariance by rototranslations of the algorithm.}, we get a smooth function defined on $\R^2$, which is generically a Morse function:
$$
f_c=\PPP_c\ast G(\sigma_x,\sigma_y).
$$
We recall that a smooth function $f_c:\R^2\to\R$ is said to be  Morse 
if it has only isolated critical points with nondegenerate Hessian. Roughly speaking, a Morse function is a function whose level sets are locally like those of Figure \ref{morse}. 

\vskip 1cm\noindent
{\bf STEP 2: The lift of $f_c:\R^2\to\R$ to a function  $\bar f_c:PT\R^2\to\R$}

This is made by associating to every point $(x,y)$ of $\R^2$ the direction $\th\in\R/(\pi\Z)$ of the level set of $f_c$ at the point $(x,y)$. This direction is well defined only at points  where $\nabla f_c\neq0$. At points where  $\nabla f_c=0$, we associate all possible directions (see Figure \ref{f-llift}).
\ffoot{nota che la funzione $f(x,y)= y^3$ ha un lift molto strano. In un certo senso qusto lift funziona bene solo per morse functions. forse bisognerebbe definire col limite quando si puo' e se nocoe sopra. comunque il vero interesse e' sulle morse functions}
More precisely, we define the {\it lifted support} $S_f$, associated with the function $f$ as follows, 
\bqn
S_f&=&\{(x,y,\th)\in \R^2\times P^1\mbox{ s.t. }~ \nabla f_c(x,y)\cdot(\cos(\th),\sin(\th))=0 \},
\eqnn
where the dot means the standard scalar product on $\R^2$. 
Let $\Pi:S_f\to \R^2$ be the standard projection $(x,y,\th)\in S_f\to (x,y)\in\R^2$. 
Notice that if $\nabla f_c(x,y)\neq0$ then $\Pi^{-1}(x,y)$ is a single point, while  if 
$\nabla f_c(x,y)=0$ then $\Pi^{-1}(x,y)=\R/(\pi\Z)$. 

\ppotR{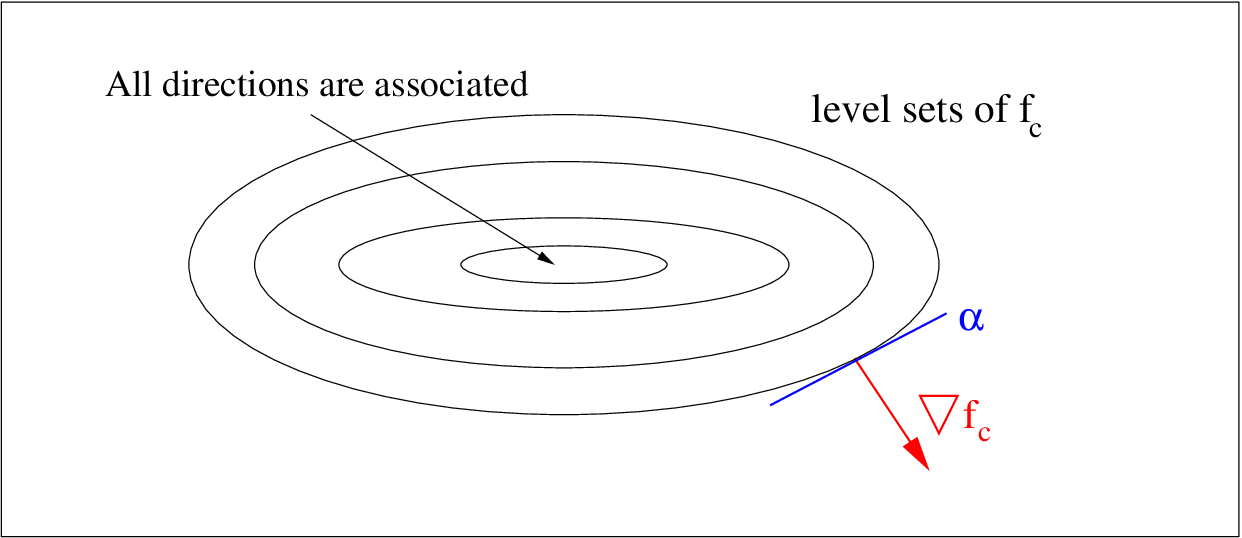}{Lift of an image with a maximum point.}{10}


Let us study the set $S_f$, when $f_c$ is a Morse function. If $(x,y)\in\R^2$ is such that $\nabla f_c(x,y)\ne0$ and $U$ is a small enough open neighborhood of $(x,y)$, then the lift of $S_f$ is an orientable manifold in $U\times P^1$. See Figure \ref{fig:f-morse-tagliata} A.  If $(x,y)$ is an isolated maximum of  $f_c$, and $U$ is a small enough open neighborhood of $(x,y)$ having a level set of $f_c$ as boundary, then $S_f$ is a M\"obius strip in $U\times P^1$. See Figure \ref{fig:f-morse-tagliata} B. The same happens when $(x,y)$ is an isolated minimum or saddle point  of  $f_c$. Indeed we have:

\immagine[12]{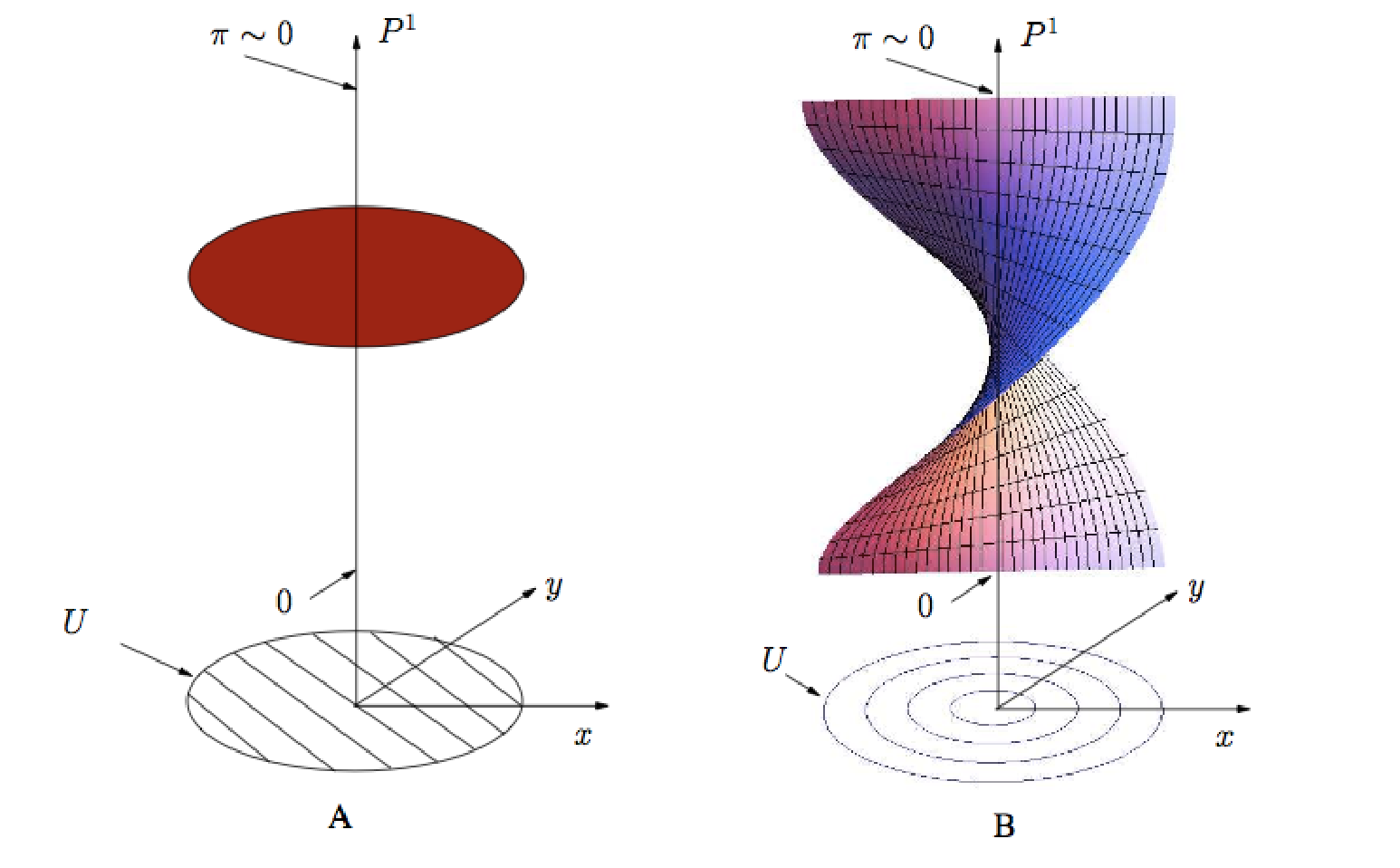}{In Figure A, we draw the support of the lift of a linear function with nonvanishing gradient. Figure $B$ presents the  support of the lift of a function in a neighborhood of a maximum point.}

\bp
\label{p-morse}
If  $f_c:\R^2\to\R$ is a Morse function, then  $S_f$ is an embedded 2-D submanifold of $\R^2\times P^1$.
\ep

\bproof
Consider the surface $\bar S\in SE(2)$ given by the equation
\bqn
g(x,y,\th):=\cos(\th)\partial_x f_c+\sin(\th)\partial_y f_c=0.
\eqn 
If $(x,y,\th)\in \bar S$ then $(x,y,\th+\pi)\in \bar S$ as well and $\bar S$ is a double covering of $S_f$. It is enough to show that $\bar S$ is a surface.

At points $(x,y,\th)\in \bar S$ where $\nabla f_c\neq0$ then $dg$ is non-zero. Indeed $\partial_\th g=0$ would imply that the vector $\nabla f_c\neq0$ is orthogonal to two non-zero vectors.  At points where  $\nabla f_c=0$ we have 
\bqn
\left(\ba{c}
\partial_x g\\
\partial_y g
\ea
\right)=H_{f_c}\cdot
\left(\ba{c}
\cos(\th)\\
\sin(\th)
\ea
\right).
\eqn
which cannot be zero since the Hessian $H_{f_c}$ of $f_c$ is non-degenerate by the Morse assumption.
\eproof

\vspace{1cm}\noindent
{\bf STEP 3: Lift of $f_c$ to a distribution in $\R^2\times P^1$ supported on $S_f$}\\\\
Consider the distribution on $\R^2\times P^1$:
$$
\bar f_c(x,y,\th)=f_c(x,y)\delta(g)
$$
where $\delta(g)$ is the Dirac-delta distribution associated with $g(x,y,\th):=\cos(\th)\partial_x f_c+\sin(\th)\partial_y f_c$ in the sense of \cite[p.222]{Gelfand}. This distribution is supported on $S_f$ and it is canonically defined by $f_c$. Notice that this choice is not crucial and there are other possibilities. \cecco{For example, in \cite{sanguinetti} the Dirac delta is replaced by a a power of the cosine of the angle, centered on the angle $\th$}.


\brem
This step is formally necessary for the following reason. The surface $S_f$ is 2D in a 3D manifold, hence the real function $f_c$ defined on it is vanishing a.e. as a function defined on $\PTR$. Thus the hypoelliptic evolution of $f_c$ (that is, the next STEP 4) produces a vanishing function.
Multiplying  $f_c$  by a Dirac delta is a natural way to obtain a nontrivial evolution.
\erem

\vspace{1cm}\noindent
{\bf STEP 4: Hypoelliptic evolution}\\\\
Fix $T>0$. Compute the solution at time $T$ to the Cauchy problem,
\bqn
\left\{\ba{l}
\partial_t \phi(x,y,\th,t)=
((\cos(\th)\partial_x+\sin(\th)\partial_y)^2+\beta^2\partial_\th^2)\phi(x,y,\th,t)\\
\phi(x,y,\th,0)=\bar f_c(x,y,\th).
\ea\right.
\eqnl{e:alg-evol}

\brem
\label{rem-beta} In the formula above, the Laplacian is given by $X_1^2+\beta^2 X_2^2$, thus it depends on the fixed parameter $\beta$. This means that we use evolution depending on the cost $J_\beta$ rather than $J$. Tuning this parameter will provide better results of the reconstruction algorithm.
\erem

\vspace{1cm}\noindent
{\bf STEP 5: Projecting down}\\\\
Compute the reconstructed image by choosing the maximum value on the fiber.
\bqn
f_T(x,y)=\max_{\th\in P^1} \phi(x,y,\th,T).
\eqnn
Again other choices are possible for this projection.\\

\brem
\label{r-reiteriagite}
The algorithm depends on two parameters. The first is the time of the evolution $T$, the second is the relative weight $\beta$ in formula \r{e:alg-evol}.
A variant of this algorithm consists of  re-iterating  the steps above  for  very short diffusion times.  This idea was already presented in \cite{citti-sarti} to  build a minimal surface.

\erem

\brem \label{r-global}
One main feature of this algorithm is that it does not need the knowledge of the corrupted part. As a consequence  the diffusion acts also in the noncorrupted region. The larger the diffusion time, the more modified image in the non-corrupted region. This is very visible comparing Figure \ref{fig:ex-1} (small diffusion time) and Figure \ref{fig:ex-2} (large diffusion time). This is one of the weak points of this completion process, that certainly takes place in the V1 cortex as a low-level process. It is the counterpart of the global character of the method.

However,  due to the highly nonisotropic character of the diffusion, this effect is not too visible from a global point of view.

Modifications of the algorithm which keep the original image unmodified are suggested in \cite{citti-sarti,sanguinetti}, by admitting the diffusion in the corrupted part only. Also, in the standard PDE-based image processing inpainting algorithms, this problem disappears. Indeed,  one solves a (stationary) elliptic-like problem on the corrupted part with Neumann boundary conditions, not an evolution equation. See e.g. \cite{PDE-inpainting}.
\erem

\section{Numerical implementation and results}
\label{s-numerical}
First we present the main lines of the algorithm used in our simulations. 

\subsection{STEPS 2-3: Lift of an image}
The formal definition of the lifted function is hard to realize numerically for two reasons: the discretization of the angle variable $\th$ and the presence of a delta function.

Both issues are solved changing the definition of the lifted function:
$$\bar{f}_c(x,y,\th)=f_c(x,y) \phi(\nabla f_c(x,y),\th),$$
where $\phi(0,\th)=1/(2\eps)\, \forall\,\th\in\R/(\pi\Z)$ and $\phi(v,\th)=\phi^1(\arg(v)-\th)$ where $\phi^1(\beta)$ is the $\pi$-periodic function assuming the following values over the interval $\Pq{0,\pi}$:
\bqn
\phi^1(\beta):=
\begin{cases}
1/(2\eps) &\mbox{ if } \beta\in\Pq{\frac{\pi}{2}-\eps,\frac{\pi}{2}+\eps},\\
0 &\mbox{otherwise,}
\end{cases}
\eqnn
for a fixed $\eps>0$.

Since the space is discretized, the non-zero values of $\bar{f}_c$ are no longer defined over a set of null measure, hence the discretized hypoelliptic diffusion gives non vanishing function for all $T>0$. Thus, it is not necessary to perform STEP 3.

\subsection{STEP 4: Hypoelliptic evolution}

\newcommand{\opic}{o\Pt{\lam^2+t}}
In this section we give the crucial ideas to compute efficiently the hypoelliptic evolution \r{e:alg-evol}.  Here $\langle.,. \rangle$ is the scalar product in $\R^2$ and $R_\th$ is the rotation operator of angle $\theta$.

 First of all, the main feature of the noncommutative Fourier transform is to desintegrate the regular representation of $SE(2)$. This was the main ingredient of the computation of the hypoelliptic heat kernel in \cite{nostro-kern}. Using the Fourier transform again, the hypoelliptic heat equation is transformed into a family of parabolic equations. These are more suitable for standard numerical methods.


Roughly speaking, the non-commutative Fourier transform $\hat f(\Lambda)$ of the function $f(x,y,\th)=:f(X,\th)$, for $\Lambda\in\R^2$, is an operator meeting:
\bqn
[\hat f(\Lambda)\psi](\th)&=&\int_{\R^2}  \int_{S^1}  f(X,\al) \psi(\al+\th) d\al   e^{2\pi i\langle R_{-\th} \Lambda,X\rangle}dX
=\widetilde{(f\ast_\th \psi )}(R_{-\th}\Lambda),
\eqn
where $\ast_\th$ is the convolution with respect to the angular variable and $\widetilde{\qquad}$ is the 2-D Fourier transform with respect to the spatial variables $X=(x,y)$.

Then it is natural to consider the Fourier transform with respect to $X$. Indeed, apply this transform $u\to\tilde{u}$ to the initial value problem: 
\bqn
\left\{\ba{l}
\partial_t u=\Delta_H u\\
u(0,X,\th)=\bar f_c(X,\th),
\ea\right.
\eqn
that gives
\bqn
\left\{\ba{l}
\partial_t \tilde u=\beta^2 \partial_\th^2\tilde{u}-4\pi^2(x \cos(\th)+y\sin(\th))^2 \tilde{ u}\\
\tilde u(0,X,\th)=\widetilde{{\bar f}_c}(X,\th).
\ea\right.
\eqnl{e-finfin}
Hence, for each point 
 in the Fourier space, we have to solve an evolution equation with Mathieu right-hand term. 
 
It is easy to solve explicitly \r{e-finfin} over $\PTR$, i.e. with $\th\in \R/\pi$. This simply divides the computation time by 4.

This is the principle of the algorithm, which is massively parallelizable, since we can solve simultaneously the equation \r{e-finfin} at each point of the Fourier space.

 \subsection{Results of image reconstruction}\label{s-risultiagite}

In this section we provide results of image reconstruction using the algorithm presented above. For these examples, we have  tuned the parameters $\beta$, that is the relative weight, and $T$, the final time of evolution.

Notice again that this algorithm processes the image globally and  does not need the information about where the image is corrupted. The counterpart is that it modifies the non-corrupted part too.

We present three results. 

\bi
\i Figure  \ref{fig:ex-1} shows an image which is corrupted in a small piece of it. Then the diffusion can be applied for a rather small time avoiding an important diffusion effect in the noncorrupted part.
 \i   Figure  \ref{fig:ex-2} shows a strongly corrupted image. In this case a larger diffusion time is necessary to ``inpaint'' completely the corrupted part.   The diffusion effect is clearly much more important. However in our opinion the result is surprisingly good.
\i The residual vertical and horizontal stripes on Figure \ref{fig:ex-2} are not due to numerical discretization (they do not occur in Figure \ref{fig:ex-1}). They are the result of the diffusion of the original (white) grid. This is again a consequence of the fact that the diffusion process is global, as explained in Remark \ref{r-global}. In the spirit of global completion, this drawback is more or less unavoidable.
 \i Due to pixelization of the image, one could think that corruption along the diagonal is the worst situation. Figure  \ref{fig:ex-3} show that this is not the case.
\ei

\immagine[11]{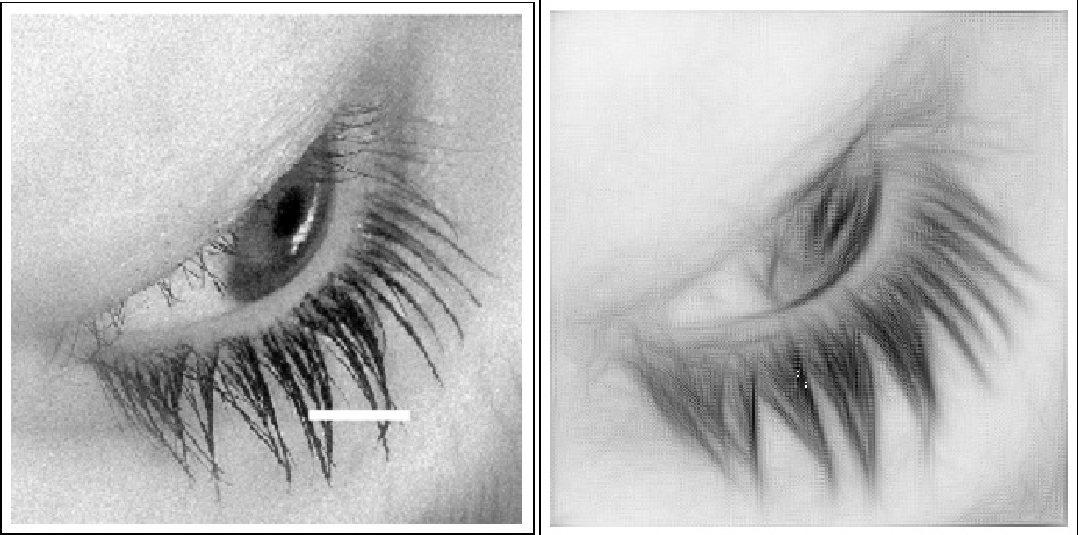}{Reconstruction of an image corrupted on a small portion. Here the diffusion is applied for a small time}
\immagine[11]{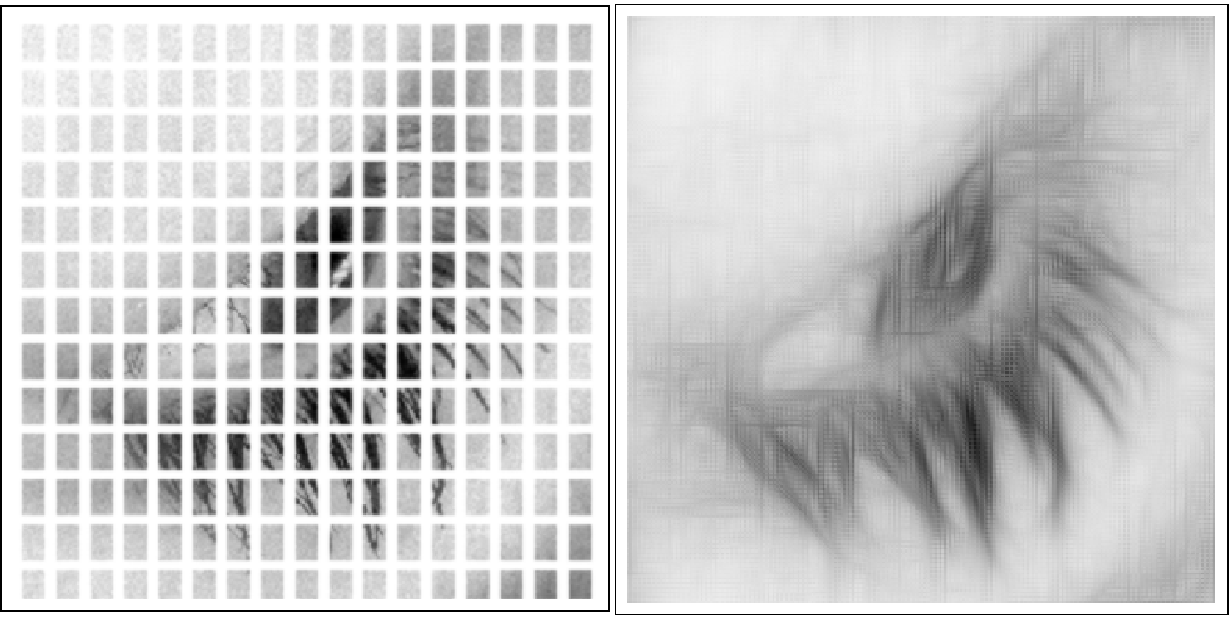}{Reconstruction of an image deeply corrupted. A larger time of diffusion is necessary}
\immagine[9]{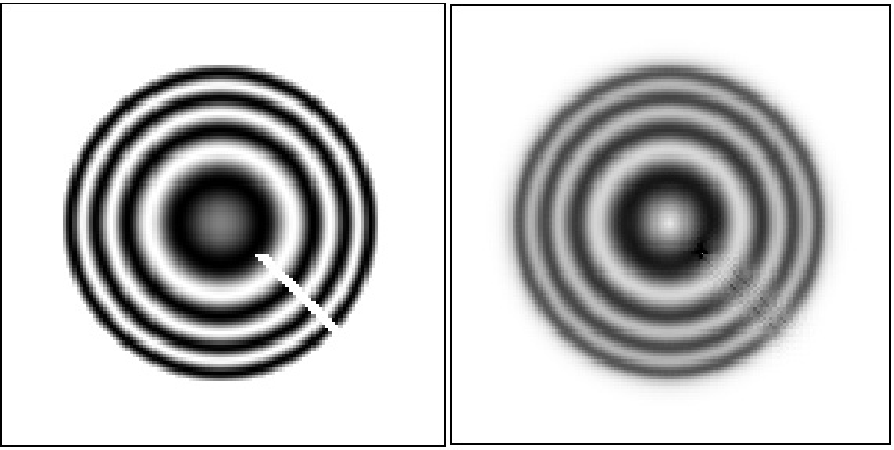}{Reconstruction of an image corrupted on the diagonal}

\section*{A Genericity of Morse properties of Gaussian convolution}

\newcommand{\myf}{\PPP}
\newcommand{\myOmega}{\OOO}

In this appendix, we prove that, generically, the convolution of a $L^2$ function over a bounded domain $\myOmega\subset\R^2$ with a Gaussian $G$ is a Morse function. In particular, we first prove in Theorem \ref{t-gen-1} that the set of functions $\myf\in L^2(\myOmega,\R)$ the convolution of which with a Gaussian is a Morse function is residual\footnote{We recall that a subset of a topological space is residual when it is a countable intersection of open and dense sets.} in $L^2(\myOmega,\R)$. We then prove in Theorem \ref{t-gen-2} that the set of functions $\myf\in L^2(\myOmega,\R)$ such that  $\myf\ast G$ {\it restricted to a compact} $K\subset\R^2$ is a Morse function is open and dense.

\newcommand{\residual}{\mathbb{X}}
\newcommand{\abres}{\mathbb{A}}
\newcommand{\eval}{\mathrm{ev}_\rho}
\newcommand{\ti}{\;\;\makebox[0pt]{$\top$}\makebox[0pt]{$\cap$}\;\;}
\newcommand{\conv}{\Gamma}

\bdeff
\label{transv}
Let $Z,Y$ be $C^1$ manifolds, $\fz{F}{Z}{Y}$ a $C^1$ map and $W\subset Y$ a submanifold. We say that $F$ is transversal to $W$ at $z\in Z$, in symbols $F\ti_z W$, if, where $y=F(z)$, either $y\not\in W$ or $y\in W$ and
\be
\i the inverse image $(T_z F)^{-1} (T_y W)$ splits and
\i the image $(T_zF)(T_zZ)$ contains a closed complement to $T_yW$ in $T_yY$.
\ee
We say that $F$ is transversal to $W$, in symbols $F\ti W$, if $F\ti_zW$ for every $z\in Z$.
\edeff

We recall that a closed subspace $F$ of a Banach space $E$ splits when there exists a closed subspace $G$ such that $E=F\oplus G$.
\brem \label{split-hilb} If $E$ is Hilbert, then every closed subspace splits. See \cite[Prop. 2.1.15]{abr-mar}.
\erem

\bt
\label{t-gen-1}
Let $\myOmega$ be a bounded domain of the plane $\R^2$. Fix $\sigma_x,\sigma_y>0$. Consider the convolution map\footnote{$\myf$ is considered to be zero outside $\OOO$.} \mmfunz{\conv}{L^2(\myOmega,\R)}{\con^\infty(\R^2)}{\myf}{\myf\ast G,} where $G$ is the  Gaussian centred at $(0,0)$
\bqn
G(x,y):=\frac1{2 \pi \sigma_x \sigma_y} e^{-\frac{x^2}{2\sigma_x^2}-\frac{y^2}{2\sigma_y^2}}.
\eqnn
Let $\residual:=\Pg{\myf\in L^2(\myOmega,\R)\mbox{ s.t. }\conv(\myf)\mbox{ is a Morse function}}$. Then, $\residual$ is residual in $L^2(\myOmega,\R)$.
\et
\bproof
The proof relies on  parametric transversality Theorems. The version we use is Abraham's formulation, see \cite[Th. 19.1]{abra}, recalled in the following.
\bt
\label{abra}
Let $\abres,X,Y$ be $\con^r$ manifolds, $\fz{\rho}{\abres}{\con^r(X,Y)}$ a $\con^r$ representation, $W\subset Y$ a submanifold, and $\fz{\eval}{X\times \abres}{Y}$ the evaluation map. Define $\abres_W\subset \abres$ by $\abres_W=\Pg{a\in \abres\ |\ \rho_a \ti W }$. Assume that:
\be
\i $X$ has a finite dimension $n$ and $W$ has a finite codimension $q$ in $Y$,
\i $\abres$ and $X$ are second countable,
\i $r>\max\Pg{0,n-q}$,
\i $\eval\ti W$.
\ee
Then, $\abres_W$ is residual in $\abres$.
\et
We apply Theorem \ref{abra} with $\abres=L^2(\myOmega,\R)$, $X=\R^2$, $r=2$. We choose $Y=\R^2\times\R\times\R^2\times\R^3$ and $\rho$ the 2-jets of $\Gamma(\myf)$, i.e.,
\mmfunz{\rho}{\abres}{\con^r(X,Y)}{\myf}{(\Pi_1,\Pi_2,\conv(\myf),\partial_x\conv(\myf),\partial_y\conv(\myf),\partial^2_{xx}\conv(\myf),\partial^2_{xy}\conv(\myf),\partial^2_{yy}\conv(\myf))}
where
\begin{center}
 \mfunz{\Pi_1}{X}{\R}{(x,y)}{x}\hspace{2mm} and \hspace{2mm} \mfunz{\Pi_2}{X}{\R}{(x,y)}{y}
\end{center}
are the canonical projections.

We fix 
\bqn
W&=&\left\{(x,y,a,p_1,p_2,q_1,q_2,q_3)\in Y\ \mbox{s.t.~~} (x,y)\in\R^2,\ p_1=p_2=0,\ q_1q_3-q_2^2=0 \right\}.
\eqnn
 A function $\myf\in C^2(\R^2)$  is a Morse function if and only if 
\bqn
\eval(x,y,\myf)=\rho_\myf\,(x,y)&=&(x,y,\Gamma(\myf)(x,y),\partial_x \Gamma(\myf)(x,y),\partial_y\Gamma(\myf)(x,y),\nn
&&\partial^2_{xx}\Gamma(\myf)(x,y),\partial^2_{xy}\Gamma(\myf)(x,y),\partial^2_{yy}\Gamma(\myf)(x,y))
\nonumber
\eqn
does not belong to $W$ for all $(x,y)\in \R^2$.

Remark that $W$ is not a manifold. However, it is an algebraic set and hence it is a finite union of manifolds. In the following, we apply Theorem \ref{abra} as if $W$ were a manifold, with the understanding that the Theorem is applied to each component.

We now verify each of the conditions 1-4 in Theorem \ref{abra}. Condition 1 holds with $n=2$ and $q\geq 3$ for each component of $W$. Condition 2 holds, since $\abres$ and $X$ are separable metric spaces and hence second countable. Condition 3 holds for each component of $W$.

Now we verify condition 4, that is the transversality condition $\eval\ti W$. Fix $x,y,\myf$ such that ${\eval}(\myf,(x,y))\in W$. Condition 1 in Definition \ref{transv} holds because of Remark \ref{split-hilb}. We now verify condition 2 in Definition \ref{transv}, where $Z=\R^2 \times \abres$. In the following, we prove that $(T_{(x,y,\myf)} \eval) (T_{(x,y,\myf)}(\R^2 \times \abres))$ is the whole $T_{{\eval}(x,y,\myf)} Y$. The map $T_{(x,y,\myf)} \eval$ has the following triangular form
\bqn
T_{(x,y,\myf)} \eval=\Pt{\ba{cc|c}
1 &0 & \ast\\
0 &1 &\ast\\
\hline
0&0& T_{\myf} \eval(x,y,\myf)
\ea}
\eqn
We are left to prove that the tangent mapping $T_{\myf} \eval(x,y,\myf)$ is surjective in $\R\times\R^2\times \R^3$, for arbitrary  $(x,y)$ fixed. After a suitable change of coordinate, we can assume that $\sigma_x=\sigma_y=1$ and that $(0,0)\in\myOmega$. Let $\eps>0$ such that $\myOmega\supset {\bf Q}:=[-\eps,\eps]\times[-\eps,\eps]$. Define the function in $L^2(\myOmega,\R)$  
$$\delta \myf(\bar x, \bar y)=\frac{c_0 +c_1 \bar x +c_2 \bar y + c_3 \bar x^2 +c_4 \bar x\bar y + c_5 \bar y^2}{G(x-\bar x,y-\bar y)}$$ restricted to ${\bf Q}$, and zero in $\myOmega\backslash {\bf Q}$. The map $\rho$ is linear in $\myf$, thus $T_\myf\eval(x,y,\myf)\Pq{\delta \myf}=\eval(x,y,\delta \myf)$. Consider the linear operator
\bqn
\eval(x,y,\delta \myf)&=&\Pt{\ba{l}\int_{\bf Q} \delta \myf(\bar x,\bar y) G(x-\bar x,y-\bar y)\, d\bar x d \bar y\\
\int_{\bf Q} \delta \myf(\bar x,\bar y) \partial_1 G(x-\bar x,y-\bar y)\, d\bar x d \bar y\\
\int_{\bf Q} \delta \myf(\bar x,\bar y) \partial_2 G(x-\bar x,y-\bar y)\, d\bar x d \bar y\\
\int_{\bf Q} \delta \myf(\bar x,\bar y) \partial^2_{11} G(x-\bar x,y-\bar y)\, d\bar x d \bar y\\
\int_{\bf Q} \delta \myf(\bar x,\bar y) \partial^2_{12} G(x-\bar x,y-\bar y)\, d\bar x d \bar y\\
\int_{\bf Q} \delta \myf(\bar x,\bar y) \partial^2_{22} G(x-\bar x,y-\bar y)\, d\bar x d \bar y\ea}
\eqnn as a function of the 6 variables $(c_0,\ldots,c_5)$, and consider the linear system $\eval(x,y,\delta \myf)=(a,p_1,p_2,q_1,q_2,q_3),$ where $(a,p_1,p_2,q_1,q_2,q_3)\in Y$ is fixed. A direct computation shows that the determinant of the system is $\frac{65536 \eps^{28}}{164025 \sigma_x^8 \sigma_y^8}>0$, thus the system always has a solution, i.e. $T_{\myf} \eval(x,y,\myf)$ is surjective.

By applying Theorem \ref{abra}, we get $\abres_W$ residual in $\abres$. We now prove that $\abres_W=\residual$. Since $\myf\in\residual$ implies $\eval(x,y,\myf)\not \in W$, then $\rho_\myf\ti W$, hence $\abres\supset\residual$.

Now let us prove the inclusion $\abres \subset \residual$. Let $\myf\in\abres$ and fix $(x,y)\in\R^2$.\\

\b{Nonintersection claim :} $\rho_\myf(x,y)\not\in W$.\\\\
{\it Proof of the claim.} By contradiction, let $$w=\eval(x,y,\myf)\in W.$$ Since $\rho_\myf\ti_{(x,y)} W$, then $\Pt{T_{(x,y)} \rho_\myf}\Pt{T_{(x,y)}\R^2}$ contains a closed complement to $T_wW$ in $T_w Y$.

Observe that $$\mathrm{dim}\Pt{T_{(x,y)} \rho_\myf}\Pt{T_{(x,y)}\R^2}\leq \mathrm{dim}\Pt{T_{(x,y)}\R^2}=2$$ and $\mathrm{codim}\,{T_wW}\geq 3$, thus $\Pt{T_{(x,y)} \rho_\myf}\Pt{T_{(x,y)}\R^2}$ cannot contain a closed complement to $T_wW$ in $T_w Y$. A contradiction.\\

By applying the claim for each $(x,y)\in \R^2$, we get that $\rho_\myf$ is a Morse function.
\eproof

\bt
\label{t-gen-2}
Let $K$ be a compact subset of $\R^2$ with non-empty interior. Under the hypothesis of Theorem \ref{t-gen-1}, the set $\residual_K:=\Pg{\myf\in L^2(\myOmega,\R) \mbox{ s.t. }{\rho_\myf}_{|_K} \mbox{ is a Morse function}}$ is open and dense in $L^2(\myOmega,\R)$.
\et
\bproof
Applying the openness of nonintersection Theorem \cite[Th. 18.1]{abra} and using the nonintersection claim, we get that $\residual_K$ is an open subset of $L^2(\myOmega,\R)$. Since $\residual_K\supset\residual$ and $\residual$ is dense, then the conclusion holds.
\eproof

\end{document}